\documentclass[a4paper,10pt]{article}
\usepackage{geometry}
\usepackage[english,french]{babel}
\usepackage{paralist}
\usepackage{pifont}
\usepackage{amsmath,amsfonts,amssymb,amsthm,cases}
\usepackage{empheq}
\usepackage{array,dcolumn,booktabs}
\usepackage[section]{placeins}
\usepackage[dvipsnames,svgnames]{xcolor}
\usepackage{tikz,todonotes}
\usepackage{subfig}
\usepackage{url}
\usepackage{gensymb}
\usepackage{microtype}
\usepackage{marvosym}
\usepackage{graphicx}
\usepackage{hyperref}
\usepackage{epstopdf}
\graphicspath{{images/}}
\usepackage{listings}

\newcommand{\ds}{\displaystyle}
\newcommand{\bg}{\begin{equation}}
\newcommand{\ed}{\end{equation}}
\usepackage{xspace}
\makeatletter
\DeclareRobustCommand{\freefem}
{\valign{\vfil\hbox{##}\vfil\cr
   \textsf{FreeFem\kern-.1em}\cr
   $\hbox{\fontsize{\sf@size}{0}\textbf{+\kern-0.05em+}}$\cr}\xspace%
}
\hypersetup{
backref=true, %permet d'ajouter des liens dans...
pagebackref=true,%...les bibliographies
hyperindex=true, %ajoute des liens dans les index.
colorlinks=true, %colorise les liens
breaklinks=true, %permet le retour \`a la ligne dans les liens trop longs
urlcolor= red, %couleur des hyperliens
linkcolor= black, %couleur des liens internes
citecolor = blue, 
 bookmarks=true, %cr\'e\'e des signets pour Acrobat
bookmarksopen=true, %si les signets Acrobat sont cr\'e\'es,
pdftitle={Propagation of a Tsunami wave}, %informations apparaissant dans
pdfauthor={Georges SADAKA}, %dans les informations du document
pdfsubject={boussinesq_freefem_article} %sous Acrobat.
} 
\makeatother
\theoremstyle{plain}

\theoremstyle{definition}

\theoremstyle{remark}
\newtheorem{rem}{Remark}

\begin{document}
{\selectlanguage{english} 
\date{}
\title{ {\sc \textbf{Generation and propagation of a {\it Tsunami} wave : a new mesh adaptation technique}}}
\author{{\sc \Large{Georges Sadaka}}\footnote{Universit\'e de Picardie Jules Verne, LAMFA CNRS UMR 7352, 33, rue Saint-Leu, 80039 Amiens, France, \url{http://lamfa.u-picardie.fr/sadaka/}  \Email\xspace \small{georges.sadaka@u-picardie.fr}}}
\maketitle
\begin{flushleft}
\textbf{Abstract} : The importance of the study of the propagation of a {\it Tsunami} wave came from the complex phenomenon and its natural disasters which represents a major risk for populations. To model this phenomena, we will consider a {\it simplified} Boussinesq\footnote{\cite{Bou1871}\label{Bou18710}} system of BBM\footnote{Benjamin, Bona and Mahony (BBM)} type (sBBM) derived by D. Mitsotakis in \cite{Mit09}, \label{Mit090} over a flat bottom then over a variable bottom in space and in time and apply this system, first, using a mesh generated using a photo of the Mediterranean sea, second, using a mesh generated using an imported xyz bathymetry for the sea near Java island and then we will consider a realistic example of the {\it Tsunami} wave near Java island which happened in 2006.\\

We choose here to use \freefem \cite{HecPioLehOht12}\label{HecPioLehOht1001} software which simplifies the construction of the domain, in particular, one of the advantage of \freefem is that we can build a mesh using a photo and we can easily export bathymetric data in order to consider more realistic simulations where a special adapt mesh technique applied for these two methods is detailed in the sequel.\\
\end{flushleft}
\tableofcontents
\section{Introduction}
We consider here the numerical simulation of the sBBM in 2D over a variable bottom in space $\textcolor{blue}{D}(x,y)$ and in time $\textcolor{red}{\zeta}(x,y,t)$ :\\
\bg\label{derivETAV}
\begin{array}{l}\eta_t+\nabla\cdot\left((\textcolor{blue}{D}+\eta+\textcolor{red}{\zeta}) V\right)+\textcolor{red}{\zeta_t}+\tilde{A}\nabla\cdot\left(\textcolor{blue}{D^2}\nabla \textcolor{red}{\zeta_t}\right)
+\nabla\cdot\left\{A\textcolor{blue}{D^2}\left[\nabla\left(\nabla \textcolor{blue}{D}\cdot V\right)+\nabla \textcolor{blue}{D}\nabla\cdot V\right]-b\textcolor{blue}{D^2}\nabla\eta_t\right\}=0,\\ \\
V_t+g\nabla\eta+\frac{1}{2}\nabla|V|^2+ Bg\textcolor{blue}{D}\left[\nabla\left(\nabla \textcolor{blue}{D}\cdot \nabla\eta\right)+\nabla \textcolor{blue}{D}\Delta\eta\right]-d\textcolor{blue}{D^2}\Delta V_t- B\textcolor{blue}{D}\nabla \textcolor{red}{\zeta_{tt}}=0,
\end{array}\ed
where 
$$\widehat{a}=\left(\theta-\frac{1}{2}\right),\qquad \widehat{b}=\frac{1}{2}\left((\theta-1)^2-\frac{1}{3}\right),\qquad\tilde{A}=\nu \widehat{a}-(1-\nu)\widehat{b}, \qquad A=-\widehat{b},\qquad B=1-\theta,$$
$$b=\frac{1}{2}\left(\theta^2-\frac{1}{3}\right)(1-\nu),\qquad d=\frac{1}{2}\left(1-\theta^2\right)(1-\mu)\qquad\mbox{ and }\qquad g=9.81 \mbox{ is the gravity}.$$
This system is an approximation to the three-dimensional Euler equations describing the irrotational free surface flow of an ideal fluid $\Omega \subset\mathbb{R}^3$, which is bounded below by $-h(x,y,t)=-\textcolor{blue}{D}(x,y)-\textcolor{red}{\zeta}(x,y,t)$ and above by the free surface elevation $\eta(x,y,t)$ (cf. Figure \ref{domomega_t}).\\
\begin{figure}[!htb]
\begin{center}
\includegraphics[height=6cm,width=8cm]{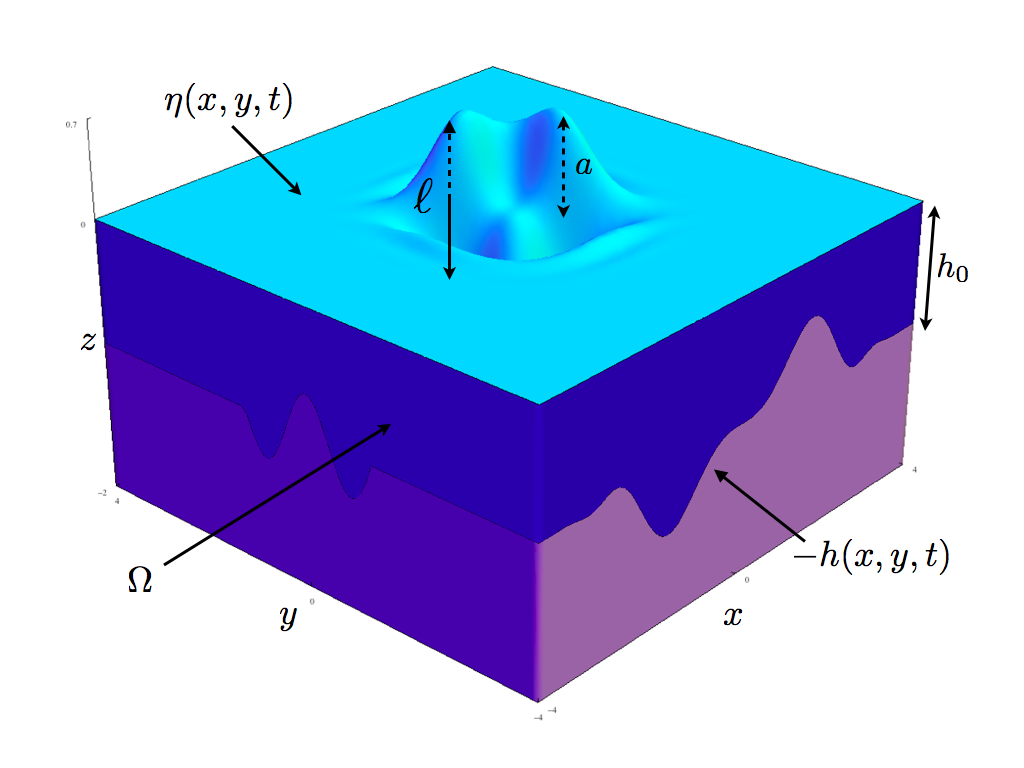}
\caption{The domain $\Omega.$}
\label{domomega_t}
\end{center}
\end{figure}
The variables in (\ref{derivETAV}) : $\textbf{X}=(x,y)\in\Omega$ and $t>0$ are proportional to position along the channel and time, respectively. $\eta=\eta(\textbf{X},t)$ being proportional to the deviation of the free surface departing from its rest position and $V=V(\textbf{X},t)=\left(\begin{array}{c} u(\textbf{X},t) \\  v(\textbf{X},t)\end{array}\right)=(u,v)^T=(u;v)$ being proportional to the horizontal velocity of the fluid at some height. $\nabla\boldsymbol{\cdot}=\left(\begin{array}{c}\partial_x\boldsymbol{\cdot} \\ \partial_y\boldsymbol{\cdot}\end{array}\right)$ is the gradient, $\nabla\cdot (\,\star\,\,;\boldsymbol{\cdot}\,)=\partial_x \star+\partial_y \boldsymbol{\cdot}$ is the divergence and $\Delta \boldsymbol{\cdot} =\partial_{xx}\boldsymbol{\cdot} +\partial_{yy}\boldsymbol{\cdot}$ is the laplacian. \\

\begin{rem}
In our study, we suppose that $\eta=\mathcal{O}(a)$, where here the amplitude $a$ is the difference between the water surface and the zero level. Also we set $\lambda=\mathcal{O}(\ell)$ be the wave length. In addition, we limit ourselves to the case where $\eta+D>0$ (there is no dry zone), since we are in a big deep water wave regime.\\
\end{rem}
This paper is organized as follows: in Section \ref{secDBS}, we present the space and time discretization of equation (\ref{derivETAV}). In Section \ref{meshgeninitdata}, we present a method to build a mesh using a photo and then using an imported $xyz$ bathymetry, in which, we will also present a special adaptive mesh technique around the initial data used for the generation of a {\it Tsunami} wave. In Section \ref{numsim}, we first check the convergence of our code, which establishes the adequacy of the chosen finite element discretization, then we simulate the propagation of a wave, that looks like a {\it Tsunami} wave generated by an Earthquake, in the Mediterranean sea over the sBBM system (\ref{derivETAV}) with a flat bottom using the mesh generated from a photo of the Mediterranean sea, then with a variable bottom in space using the mesh generated from the xyz bathymetry of the sea near Java island and finally, with a realistic example of the {\it Tsunami} wave near Java island which happened in 2006.\\
\section{Discretization of the sBBM system}\label{secDBS}
In this section, we present the spatial discretization of (\ref{derivETAV}) using finite element method with $\mathbb{P}_1$ continuous piecewise linear functions and for the time marching scheme an explicit second order Runge-Kutta scheme.
\subsection{Spatial discretization}
We let $\Omega$ be a convex, plane domain, and $\mathcal{T}_h$ be a regular, quasiuniform triangulation of $\Omega$ with triangles of maximum size $h<1$. Setting $V_h=\{v_h\in C^0(\bar{\Omega}); v_h|_T \in \mathbb{P}_1(T),\forall T\in\mathcal{T}_h\}$ be a finite-dimensional subspace of $H^1(\Omega)$, where $\mathbb{P}_1$ is the set of all polynomials of degree $\leq 1$ with real coefficients and denoting by $\left\langle\cdot ; \cdot\right\rangle$ the $L^2$ inner product on $\Omega$, we consider the weak formulation of system (\ref{derivETAV}) : 

Find $\eta_h, u_h, v_h\in V_h$ such that $\forall \phi^\eta_h,\phi^u_h,\phi^v_h\in V_h$, we have :\\
\bg\label{derivETAVweak}\begin{array}{l}
\left\langle\eta_{ht}-b\nabla\cdot\left(\textcolor{blue}{D^2}\nabla\eta_{ht}\right)+\nabla\cdot\left((\textcolor{blue}{D}+\eta_h+\textcolor{red}{\zeta}) \left(u_h;v_h\right)\right)+\textcolor{red}{\zeta_t};\phi^\eta_h\right\rangle+\left\langle\tilde{A}\nabla\cdot\left(\textcolor{blue}{D^2}\nabla \textcolor{red}{\zeta_t}\right);\phi^\eta_h\right\rangle\\
\qquad\qquad\qquad\qquad\qquad\qquad\qquad\qquad+\left\langle\nabla\cdot\left\{A\textcolor{blue}{D^2}\left[\nabla\left(\nabla \textcolor{blue}{D}\cdot \left(u_h;v_h\right)\right)+\nabla \textcolor{blue}{D}\nabla\cdot \left(u_h;v_h\right)\right]\right\};\phi^\eta_h\right\rangle=0,\\
\left\langle u_{ht}-d\textcolor{blue}{D^2}\Delta u_{ht}+g\eta_{xh}+u_h u_{hx}+ v_h v_{hx}- B\textcolor{blue}{D}\textcolor{red}{\zeta_{xtt}};\phi^u_h\right\rangle+ Bg\left\langle\textcolor{blue}{D}\Big[\left(\nabla \textcolor{blue}{D}\cdot \nabla\eta_h\right)_x+\textcolor{blue}{D}_x\Delta\eta_h\Big];\phi^u_h\right\rangle=0,\\
\left\langle v_{ht}-d\textcolor{blue}{D^2}\Delta v_{ht}+g\eta_{yh}+u_h u_{hy}+ v_h v_{hy}- B\textcolor{blue}{D}\textcolor{red}{\zeta_{ytt}};\phi^u_h\right\rangle+ Bg\left\langle\textcolor{blue}{D}\left[\left(\nabla \textcolor{blue}{D}\cdot \nabla\eta_h\right)_y+\textcolor{blue}{D}_y\Delta\eta_h\right];\phi^v_h\right\rangle=0.
\end{array}\ed

For simplicity, we set $\phi^\eta_h=\Phi^\eta,\ \phi^u_h=\Phi^u,\  \phi^v_h=\Phi^v,\ \eta_h=\mathcal{E},\ u_h=\mathcal{U},v_h=\mathcal{V}$, so that system (\ref{derivETAVweak}) can be rewrite in the following way : 
\bg\label{derivBBMFull}\hspace{-1cm}\left\{\begin{array}{rcl} \Big\langle\partial_t\mathcal{E}-b\nabla\cdot\left(\textcolor{blue}{D^2}\nabla\partial_t\mathcal{E}\right);\Phi^\eta\Big\rangle&=&-\left\langle (\textcolor{blue}{D}+\textcolor{red}{\zeta}+\mathcal{E})\nabla\cdot(\mathcal{U};\mathcal{V})+(\textcolor{blue}{D_x}+\textcolor{red}{\zeta_x}+\mathcal{E}_x)\mathcal{U}+(\textcolor{blue}{D_y}+\textcolor{red}{\zeta_y}+\mathcal{E}_y)\mathcal{V}+\textcolor{red}{\zeta_t}\right.\\
&&\left.+\tilde{A}\nabla\cdot\left(\textcolor{blue}{D^2}\nabla \textcolor{red}{\zeta_t}\right)+A\nabla\cdot\left\{\textcolor{blue}{D^2}\left[\nabla\left(\nabla \textcolor{blue}{D}\cdot (\mathcal{U};\mathcal{V})\right)+\nabla \textcolor{blue}{D}\nabla\cdot (\mathcal{U};\mathcal{V})\right]\right\};\Phi^\eta\right\rangle\\
&=&\mathcal{F}\left(\mathcal{E},\mathcal{U},\mathcal{V},\Phi^\eta\right),\\
\Big\langle(I_d-d\textcolor{blue}{D^2}\Delta)\partial_t\mathcal{U};\Phi^u\Big\rangle&=&-\left\langle g\mathcal{E}_{x}+\mathcal{U}\mathcal{U}_{x}+\mathcal{V}\mathcal{V}_{x}+Bg\textcolor{blue}{D}\Big[\left(\nabla \textcolor{blue}{D}\cdot \nabla\mathcal{E}\right)_x+\textcolor{blue}{D}_x\Delta\mathcal{E}\Big]-B\textcolor{blue}{D} \textcolor{red}{\zeta_{xtt}};\Phi^u\right\rangle\\
&=&\mathcal{G}\left(\mathcal{E},\mathcal{U},\mathcal{V},\Phi^u\right),\\
\Big\langle(I_d-d\textcolor{blue}{D^2}\Delta)\partial_t\mathcal{V};\Phi^v\Big\rangle&=&-\left\langle g\mathcal{E}_{y}+\mathcal{U}\mathcal{U}_{y}+\mathcal{V}\mathcal{V}_{y}+Bg\textcolor{blue}{D}\Big[\left(\nabla \textcolor{blue}{D}\cdot \nabla\mathcal{E}\right)_y+\textcolor{blue}{D}_y\Delta\mathcal{E}\Big]-B\textcolor{blue}{D} \textcolor{red}{\zeta_{ytt}};\Phi^v\right\rangle\\
&=&\mathcal{H}\left(\mathcal{E},\mathcal{U},\mathcal{V},\Phi^u\right).\\
\end{array}\right.\ed
Next, we discretize system (\ref{derivBBMFull}). First, integrating by parts, the left hand side in (\ref{derivBBMFull}) gives :\\
$$
-\Big\langle b\nabla\cdot\left(\textcolor{blue}{D^2}\nabla\partial_t\mathcal{E}\right);\Phi^\eta\Big\rangle=b\left\langle \textcolor{blue}{D^2}\nabla\partial_t\mathcal{E};\nabla(\Phi^\eta)\right\rangle-\ds\int_{\Gamma_n}b\textcolor{blue}{D^2}\Phi^\eta\ds\frac{\partial (\partial_t\mathcal{E})}{\partial n}\partial \gamma,
$$
$$
-\Big\langle d\textcolor{blue}{D^2}\Delta\partial_t\mathcal{U};\Phi^u\Big\rangle=d\left\langle\textcolor{blue}{D^2}\nabla\partial_t\mathcal{U};\nabla\Phi^u\right\rangle+d\left\langle2\textcolor{blue}{D}\nabla\textcolor{blue}{D}\cdot\nabla\partial_t\mathcal{U};\Phi^u\right\rangle-\ds\int_{\Gamma_n}d\textcolor{blue}{D^2}\Phi^u\ds\frac{\partial (\partial_t\mathcal{U})}{\partial n}\partial \gamma,
$$
and
$$
-\Big\langle d\textcolor{blue}{D^2}\Delta\partial_t\mathcal{V};\Phi^v\Big\rangle=d\left\langle\textcolor{blue}{D^2}\nabla\partial_t\mathcal{V};\nabla\Phi^v\right\rangle+d\left\langle2\textcolor{blue}{D}\nabla\textcolor{blue}{D}\cdot\nabla\partial_t\mathcal{V};\Phi^v\right\rangle-\ds\int_{\Gamma_n}d\textcolor{blue}{D^2}\Phi^v\ds\frac{\partial (\partial_t\mathcal{V})}{\partial n}\partial \gamma.
$$
Now, dealing with the right-hand side of the first equation in system (\ref{derivBBMFull}), we have :
$$
\left\langle \nabla\cdot\left(\textcolor{blue}{D^2}\nabla \textcolor{red}{\zeta_t}\right);\Phi^\eta\right\rangle=\left\langle \left(\textcolor{blue}{D^2}\textcolor{red}{\zeta_{xt}}\right)_x+\left(\textcolor{blue}{D^2}\textcolor{red}{\zeta_{yt}}\right)_y;\Phi^\eta\right\rangle=\left\langle 2\textcolor{blue}{D}\textcolor{blue}{D_x}\textcolor{red}{\zeta_{xt}}+\textcolor{blue}{D^2}\textcolor{red}{\zeta_{xxt}}+2\textcolor{blue}{D}\textcolor{blue}{D_y}\textcolor{red}{\zeta_{yt}}+\textcolor{blue}{D^2}\textcolor{red}{\zeta_{yyt}};\Phi^\eta\right\rangle,
$$
and
$$
\left\langle\nabla\cdot\left\{\textcolor{blue}{D^2}\left[\nabla\left(\nabla \textcolor{blue}{D}\cdot (\mathcal{U};\mathcal{V})\right)+\nabla \textcolor{blue}{D}\nabla\cdot (\mathcal{U};\mathcal{V})\right]\right\};\Phi^\eta\right\rangle
$$
$$
=\left\langle \nabla\cdot\left\{\textcolor{blue}{D^2}\left[\left(\left(\textcolor{blue}{D_x}\mathcal{U}+\textcolor{blue}{D_y}\mathcal{V}\right)_x;\left(\textcolor{blue}{D_x}\mathcal{U}+\textcolor{blue}{D_y}\mathcal{V}\right)_y\right)+\left(\textcolor{blue}{D_x}\nabla \cdot (\mathcal{U};\mathcal{V});\textcolor{blue}{D_y}\nabla \cdot (\mathcal{U};\mathcal{V})\right)\right]\right\};\Phi^\eta\right\rangle
$$
$$
=\left\langle \nabla\cdot\left(\textcolor{blue}{D^2}\textcolor{blue}{D_{xx}}\mathcal{U}+\textcolor{blue}{D^2}\textcolor{blue}{D_{x}}\mathcal{U}_{x}+\textcolor{blue}{D^2}\textcolor{blue}{D_{xy}}\mathcal{V}+\textcolor{blue}{D^2}\textcolor{blue}{D_{y}}\mathcal{V}_{x}+\textcolor{blue}{D^2}\textcolor{blue}{D_x}\nabla \cdot (\mathcal{U};\mathcal{V});\textcolor{blue}{D^2}\textcolor{blue}{D_{xy}}\mathcal{U}+\textcolor{blue}{D^2}\textcolor{blue}{D_x}\mathcal{U}_{y}\right.\right.
$$
$$
\left.\left.+\textcolor{blue}{D^2}\textcolor{blue}{D_{yy}}\mathcal{V}+\textcolor{blue}{D^2}\textcolor{blue}{D_y}\mathcal{V}_{y}+\textcolor{blue}{D^2}\textcolor{blue}{D_y}\nabla \cdot (\mathcal{U};\mathcal{V})\right);\Phi^\eta\right\rangle
$$
$$
=\left\langle (2\textcolor{blue}{D}\textcolor{blue}{D_x}\textcolor{blue}{D_{xx}}+2\textcolor{blue}{D}\textcolor{blue}{D_y}\textcolor{blue}{D_{xy}}+\textcolor{blue}{D^2}\textcolor{blue}{D_{xyy}}+\textcolor{blue}{D^2}\textcolor{blue}{D_{xxx}})\mathcal{U}+(2\textcolor{blue}{D}\textcolor{blue}{D_x}\textcolor{blue}{D_{xy}}+2\textcolor{blue}{D}\textcolor{blue}{D_y}\textcolor{blue}{D_{yy}}+\textcolor{blue}{D^2}\textcolor{blue}{D_{yyy}}\right.
$$
$$
+\textcolor{blue}{D^2}\textcolor{blue}{D_{xxy}})\mathcal{V}+(4\textcolor{blue}{D}\textcolor{blue}{D_{x}^2}+3\textcolor{blue}{D^2}\textcolor{blue}{D_{xx}}+2\textcolor{blue}{D}\textcolor{blue}{D_y^2}+\textcolor{blue}{D^2}\textcolor{blue}{D_{yy}})\mathcal{U}_{x}+2(\textcolor{blue}{D^2}\textcolor{blue}{D_{xy}}+\textcolor{blue}{D}\textcolor{blue}{D_x}\textcolor{blue}{D_y})\mathcal{U}_{y}+(4\textcolor{blue}{D}\textcolor{blue}{D_y^2}
$$
$$
\left.+3\textcolor{blue}{D^2}\textcolor{blue}{D_{yy}}+2\textcolor{blue}{D}\textcolor{blue}{D_x^2}+\textcolor{blue}{D^2}\textcolor{blue}{D_{xx}})\mathcal{V}_{y}+2(\textcolor{blue}{D}\textcolor{blue}{D_x}\textcolor{blue}{D_{y}}+\textcolor{blue}{D^2}\textcolor{blue}{D_{xy}})\mathcal{V}_{x};\Phi^\eta\right\rangle+\left(\left\langle 2\textcolor{blue}{D^2}\textcolor{blue}{D_{x}}\mathcal{U}_{xx};\Phi^\eta\right\rangle\right.
$$
$$
\left.+\left\langle\textcolor{blue}{D^2}\textcolor{blue}{D_{y}}\mathcal{U}_{xy};\Phi^\eta\right\rangle+\left\langle\textcolor{blue}{D^2}\textcolor{blue}{D_x}\mathcal{U}_{yy};\Phi^\eta\right\rangle+\left\langle\textcolor{blue}{D^2}\textcolor{blue}{D_{y}}\mathcal{V}_{xx};\Phi^\eta\right\rangle+\left\langle\textcolor{blue}{D^2}\textcolor{blue}{D_x}\mathcal{V}_{xy};\Phi^\eta\right\rangle+\left\langle2\textcolor{blue}{D^2}\textcolor{blue}{D_y}\mathcal{V}_{yy};\Phi^\eta\right\rangle\right);
$$
On the other hand, we have :
$$
\left\langle 2\textcolor{blue}{D^2}\textcolor{blue}{D_{x}}\mathcal{U}_{xx};\Phi^\eta\right\rangle=-\left\langle 2\textcolor{blue}{D^2}\textcolor{blue}{D_{x}}\mathcal{U}_{x};\Phi^{\eta}_x\right\rangle-\left\langle (4\textcolor{blue}{D}\textcolor{blue}{D_{x}^2}+2\textcolor{blue}{D^2}\textcolor{blue}{D_{xx}})\mathcal{U}_{x};\Phi^\eta\right\rangle+\int_{\Gamma_n}2\textcolor{blue}{D^2}\textcolor{blue}{D_{x}}\Phi^\eta\ds\frac{\partial \mathcal{U}}{\partial n}\partial \gamma,
$$
$$
\left\langle\textcolor{blue}{D^2}\textcolor{blue}{D_{y}}\mathcal{U}_{xy};\Phi^\eta\right\rangle=-\left\langle \textcolor{blue}{D^2}\textcolor{blue}{D_{y}}\mathcal{U}_{x};\Phi^{\eta}_y\right\rangle-\left\langle (2\textcolor{blue}{D}\textcolor{blue}{D_{y}^2}+\textcolor{blue}{D^2}\textcolor{blue}{D_{yy}})\mathcal{U}_{x};\Phi^\eta\right\rangle+\int_{\Gamma_n}\textcolor{blue}{D^2}\textcolor{blue}{D_{y}}\Phi^\eta\ds\frac{\partial \mathcal{U}}{\partial n}\partial \gamma,
$$
$$
\left\langle\textcolor{blue}{D^2}\textcolor{blue}{D_x}\mathcal{U}_{yy};\Phi^\eta\right\rangle=-\left\langle \textcolor{blue}{D^2}\textcolor{blue}{D_{x}}\mathcal{U}_{y};\Phi^{\eta}_y\right\rangle-\left\langle (2\textcolor{blue}{D}\textcolor{blue}{D_{x}}\textcolor{blue}{D_{y}}+\textcolor{blue}{D^2}\textcolor{blue}{D_{xy}})\mathcal{U}_{y};\Phi^\eta\right\rangle+\int_{\Gamma_n}\textcolor{blue}{D^2}\textcolor{blue}{D_{x}}\Phi^\eta\ds\frac{\partial \mathcal{U}}{\partial n}\partial \gamma,
$$
$$
\left\langle\textcolor{blue}{D^2}\textcolor{blue}{D_{y}}\mathcal{V}_{xx};\Phi^\eta\right\rangle=-\left\langle \textcolor{blue}{D^2}\textcolor{blue}{D_{y}}\mathcal{V}_{x};\Phi^{\eta}_x\right\rangle-\left\langle (2\textcolor{blue}{D}\textcolor{blue}{D_{x}}\textcolor{blue}{D_{y}}+\textcolor{blue}{D^2}\textcolor{blue}{D_{xy}})\mathcal{V}_{x};\Phi^\eta\right\rangle+\int_{\Gamma_n}\textcolor{blue}{D^2}\textcolor{blue}{D_{y}}\Phi^\eta\ds\frac{\partial \mathcal{V}}{\partial n}\partial \gamma,
$$
$$
\left\langle\textcolor{blue}{D^2}\textcolor{blue}{D_x}\mathcal{V}_{xy};\Phi^\eta\right\rangle=-\left\langle \textcolor{blue}{D^2}\textcolor{blue}{D_{x}}\mathcal{V}_{x};\Phi^{\eta}_y\right\rangle-\left\langle (2\textcolor{blue}{D}\textcolor{blue}{D_x}\textcolor{blue}{D_y}+\textcolor{blue}{D^2}\textcolor{blue}{D_{xy}})\mathcal{V}_{x};\Phi^\eta\right\rangle+\int_{\Gamma_n}\textcolor{blue}{D^2}\textcolor{blue}{D_{x}}\Phi^\eta\ds\frac{\partial \mathcal{V}}{\partial n}\partial \gamma,
$$
$$
\left\langle2\textcolor{blue}{D^2}\textcolor{blue}{D_y}\mathcal{V}_{yy};\Phi^\eta\right\rangle=-\left\langle 2\textcolor{blue}{D^2}\textcolor{blue}{D_{y}}\mathcal{V}_{y};\Phi^{\eta}_y\right\rangle-\left\langle (4\textcolor{blue}{D}\textcolor{blue}{D_{y}^2}+2\textcolor{blue}{D^2}\textcolor{blue}{D_{yy}})\mathcal{V}_{y};\Phi^\eta\right\rangle+\int_{\Gamma_n}2\textcolor{blue}{D^2}\textcolor{blue}{D_{y}}\Phi^\eta\ds\frac{\partial \mathcal{V}}{\partial n}\partial \gamma,
$$
and, consequently,
$$
\mathcal{F}\left(\mathcal{E},\mathcal{U},\mathcal{V},\Phi^\eta\right)=-\left\langle (\textcolor{blue}{D}+\textcolor{red}{\zeta}+\mathcal{E})\nabla\cdot(\mathcal{U};\mathcal{V})+(\textcolor{blue}{D_x}+\textcolor{red}{\zeta_x}+\mathcal{E}_x)\mathcal{U}+(\textcolor{blue}{D_y}+\textcolor{red}{\zeta_y}+\mathcal{E}_y)\mathcal{V}+\textcolor{red}{\zeta_t};\Phi^\eta\right\rangle
$$
$$
-\tilde{A}\left\langle 2\textcolor{blue}{D}\textcolor{blue}{D_x}\textcolor{red}{\zeta_{xt}}+\textcolor{blue}{D^2}\textcolor{red}{\zeta_{xxt}}+2\textcolor{blue}{D}\textcolor{blue}{D_y}\textcolor{red}{\zeta_{yt}}+\textcolor{blue}{D^2}\textcolor{red}{\zeta_{yyt}};\Phi^\eta\right\rangle-A\left\langle (2\textcolor{blue}{D}\textcolor{blue}{D_x}\textcolor{blue}{D_{xx}}+2\textcolor{blue}{D}\textcolor{blue}{D_y}\textcolor{blue}{D_{xy}}+\textcolor{blue}{D^2}\textcolor{blue}{D_{xyy}}\right.
$$
$$
+\textcolor{blue}{D^2}\textcolor{blue}{D_{xxx}})\mathcal{U}+(2\textcolor{blue}{D}\textcolor{blue}{D_x}\textcolor{blue}{D_{xy}}+2\textcolor{blue}{D}\textcolor{blue}{D_y}\textcolor{blue}{D_{yy}}+\textcolor{blue}{D^2}\textcolor{blue}{D_{yyy}}+\textcolor{blue}{D^2}\textcolor{blue}{D_{xxy}})\mathcal{V}+\textcolor{blue}{D^2}\textcolor{blue}{D_{xx}}\mathcal{U}_{x}+\textcolor{blue}{D^2}\textcolor{blue}{D_{xy}}\mathcal{U}_{y}-2\textcolor{blue}{D}\textcolor{blue}{D_x}\textcolor{blue}{D_{y}}\mathcal{V}_{x}
$$
$$
\left.+(\textcolor{blue}{D^2}\textcolor{blue}{D_{yy}}+2\textcolor{blue}{D}\textcolor{blue}{D_x^2}+\textcolor{blue}{D^2}\textcolor{blue}{D_{xx}})\mathcal{V}_{y};\Phi^\eta\right\rangle+A\left(\left\langle 2\textcolor{blue}{D^2}\textcolor{blue}{D_{x}}\mathcal{U}_{x}+\textcolor{blue}{D^2}\textcolor{blue}{D_{y}}\mathcal{V}_{x};\Phi^{\eta}_x\right\rangle+\left\langle \textcolor{blue}{D^2}\textcolor{blue}{D_{y}}\mathcal{U}_{x}+\textcolor{blue}{D^2}\textcolor{blue}{D_{x}}\mathcal{U}_{y}\right.\right.
$$
$$
\left.\left.+\textcolor{blue}{D^2}\textcolor{blue}{D_{x}}\mathcal{V}_{x}+2\textcolor{blue}{D^2}\textcolor{blue}{D_{y}}\mathcal{V}_{y};\Phi^{\eta}_y\right\rangle\right)-A\int_{\Gamma_n}\left((3\textcolor{blue}{D^2}\textcolor{blue}{D_{x}}+\textcolor{blue}{D^2}\textcolor{blue}{D_{y}})\Phi^\eta\ds\frac{\partial \mathcal{U}}{\partial n}+(\textcolor{blue}{D^2}\textcolor{blue}{D_{x}}+3\textcolor{blue}{D^2}\textcolor{blue}{D_{y}})\Phi^\eta\ds\frac{\partial \mathcal{V}}{\partial n}\right)\partial \gamma.
$$
For the right-hand side of second equation in system (\ref{derivBBMFull}), we have :
$$
\mathcal{G}\left(\mathcal{E},\mathcal{U},\mathcal{V},\Phi^u\right)=-\left\langle g\mathcal{E}_{x}+\mathcal{U}\mathcal{U}_{x}+\mathcal{V}\mathcal{V}_{x}+Bg\textcolor{blue}{D}\left[\left( \textcolor{blue}{D_x}\mathcal{E}_x+\textcolor{blue}{D_y}\mathcal{E}_y\right)_x+ \textcolor{blue}{D_x} (\mathcal{E}_{xx}+\mathcal{E}_{yy})\right]-B\textcolor{blue}{D} \textcolor{red}{\zeta_{xtt}};\Phi^u\right\rangle
$$
$$
=-\left\langle g\mathcal{E}_{x}+\mathcal{U}\mathcal{U}_{x}+\mathcal{V}\mathcal{V}_{x}+Bg\left(\textcolor{blue}{D}\textcolor{blue}{D_{xx}}\mathcal{E}_x+\textcolor{blue}{D}\textcolor{blue}{D_{xy}}\mathcal{E}_y\right)-B\textcolor{blue}{D} \textcolor{red}{\zeta_{xtt}};\Phi^u\right\rangle-Bg\left\langle 2\textcolor{blue}{D}\textcolor{blue}{D_x}\mathcal{E}_{xx}+\textcolor{blue}{D}\textcolor{blue}{D_y}\mathcal{E}_{xy}+\textcolor{blue}{D}\textcolor{blue}{D_x}\mathcal{E}_{yy};\Phi^u\right\rangle
$$
$$
=-\left\langle g\mathcal{E}_{x}+\mathcal{U}\mathcal{U}_{x}+\mathcal{V}\mathcal{V}_{x}+Bg\left(\textcolor{blue}{D}\textcolor{blue}{D_{xx}}\mathcal{E}_x+\textcolor{blue}{D}\textcolor{blue}{D_{xy}}\mathcal{E}_y\right)-B\textcolor{blue}{D} \textcolor{red}{\zeta_{xtt}};\Phi^u\right\rangle+Bg\left\langle 2\textcolor{blue}{D}\textcolor{blue}{D_x}\mathcal{E}_{x};\Phi^u_x\right\rangle+Bg\left\langle (2\textcolor{blue}{D_x^2}+2\textcolor{blue}{D}\textcolor{blue}{D_{xx}})\mathcal{E}_{x};\Phi^u\right\rangle
$$
$$
+Bg\left\langle \textcolor{blue}{D}\textcolor{blue}{D_y}\mathcal{E}_{x}+\textcolor{blue}{D}\textcolor{blue}{D_x}\mathcal{E}_{y};\Phi^u_y\right\rangle+Bg\left\langle (\textcolor{blue}{D_y^2}+\textcolor{blue}{D}\textcolor{blue}{D_{yy}})\mathcal{E}_{x}+(\textcolor{blue}{D_x}\textcolor{blue}{D_y}+\textcolor{blue}{D}\textcolor{blue}{D_{xy}})\mathcal{E}_{y};\Phi^u\right\rangle-\int_{\Gamma_n}Bg(3\textcolor{blue}{D}\textcolor{blue}{D_x}+\textcolor{blue}{D}\textcolor{blue}{D_y})\Phi^u\ds\frac{\partial \mathcal{E}}{\partial n}\partial \gamma
$$
$$
=-\left\langle g\left(I_d-B\left(\textcolor{blue}{D}\textcolor{blue}{D_{xx}}+2\textcolor{blue}{D_x^2}+\textcolor{blue}{D}\textcolor{blue}{D_{yy}}+\textcolor{blue}{D_y^2}\right)\right)\mathcal{E}_{x}+\mathcal{U}\mathcal{U}_{x}+\mathcal{V}\mathcal{V}_{x}-Bg\textcolor{blue}{D_x}\textcolor{blue}{D_{y}}\mathcal{E}_y-B\textcolor{blue}{D} \textcolor{red}{\zeta_{xtt}};\Phi^u\right\rangle+Bg\left\langle 2\textcolor{blue}{D}\textcolor{blue}{D_x}\mathcal{E}_{x};\Phi^u_x\right\rangle\qquad
$$
$$
+Bg\left\langle \textcolor{blue}{D}\textcolor{blue}{D_y}\mathcal{E}_{x}+\textcolor{blue}{D}\textcolor{blue}{D_x}\mathcal{E}_{y};\Phi^u_y\right\rangle-\int_{\Gamma_n}Bg(3\textcolor{blue}{D}\textcolor{blue}{D_x}+\textcolor{blue}{D}\textcolor{blue}{D_y})\Phi^u\ds\frac{\partial \mathcal{E}}{\partial n}\partial \gamma.
$$
Finally, for the right-hand side of the third equation in system (\ref{derivBBMFull}), we have :
$$
\mathcal{H}\left(\mathcal{E},\mathcal{U},\mathcal{V},\Phi^v\right)=-\left\langle g\mathcal{E}_{y}+\mathcal{U}\mathcal{U}_{y}+\mathcal{V}\mathcal{V}_{y}+Bg\textcolor{blue}{D}\left[\left( \textcolor{blue}{D_x}\mathcal{E}_x+\textcolor{blue}{D_y}\mathcal{E}_y\right)_y+ \textcolor{blue}{D_y} (\mathcal{E}_{xx}+\mathcal{E}_{yy})\right]-B\textcolor{blue}{D} \textcolor{red}{\zeta_{ytt}};\Phi^v\right\rangle
$$
$$
=-\left\langle g\mathcal{E}_{y}+\mathcal{U}\mathcal{U}_{y}+\mathcal{V}\mathcal{V}_{y}+Bg\left(\textcolor{blue}{D}\textcolor{blue}{D_{xy}}\mathcal{E}_x+\textcolor{blue}{D}\textcolor{blue}{D_{yy}}\mathcal{E}_y\right)-B\textcolor{blue}{D} \textcolor{red}{\zeta_{ytt}};\Phi^v\right\rangle-Bg\left\langle \textcolor{blue}{D}\textcolor{blue}{D_y}\mathcal{E}_{xx}+\textcolor{blue}{D}\textcolor{blue}{D_x}\mathcal{E}_{xy}+2\textcolor{blue}{D}\textcolor{blue}{D_y}\mathcal{E}_{yy});\Phi^v\right\rangle
$$
$$
=-\left\langle g\mathcal{E}_{y}+\mathcal{U}\mathcal{U}_{y}+\mathcal{V}\mathcal{V}_{y}+Bg\left(\textcolor{blue}{D}\textcolor{blue}{D_{xy}}\mathcal{E}_x+\textcolor{blue}{D}\textcolor{blue}{D_{yy}}\mathcal{E}_y\right)-B\textcolor{blue}{D} \textcolor{red}{\zeta_{ytt}};\Phi^v\right\rangle+Bg\left\langle \textcolor{blue}{D}\textcolor{blue}{D_y}\mathcal{E}_{x};\Phi^v_x\right\rangle+Bg\left\langle (\textcolor{blue}{D_{x}}\textcolor{blue}{D_{y}}+\textcolor{blue}{D}\textcolor{blue}{D_{xy}})\mathcal{E}_{x};\Phi^v\right\rangle
$$
$$
+Bg\left\langle \textcolor{blue}{D}\textcolor{blue}{D_x}\mathcal{E}_{x}+2\textcolor{blue}{D}\textcolor{blue}{D_y}\mathcal{E}_{y};\Phi^v_y\right\rangle+Bg\left\langle (\textcolor{blue}{D_x}\textcolor{blue}{D_y}+\textcolor{blue}{D}\textcolor{blue}{D_{xy}})\mathcal{E}_{x}+(2\textcolor{blue}{D_y^2}+2\textcolor{blue}{D}\textcolor{blue}{D_{yy}})\mathcal{E}_{y};\Phi^v\right\rangle-\int_{\Gamma_n}Bg(\textcolor{blue}{D}\textcolor{blue}{D_x}+3\textcolor{blue}{D}\textcolor{blue}{D_y})\Phi^v\ds\frac{\partial \mathcal{E}}{\partial n}\partial \gamma
$$
$$
=-\left\langle -Bg(2\textcolor{blue}{D_x}\textcolor{blue}{D_{y}}+\textcolor{blue}{D}\textcolor{blue}{D_{xy}})\mathcal{E}_x+\mathcal{U}\mathcal{U}_{y}+\mathcal{V}\mathcal{V}_{y}+g\left(I_d-B\left(\textcolor{blue}{D}\textcolor{blue}{D_{yy}}-2\textcolor{blue}{D_y^2}\right)\right)\mathcal{E}_{y}-B\textcolor{blue}{D} \textcolor{red}{\zeta_{ytt}};\Phi^v\right\rangle+Bg\left\langle \textcolor{blue}{D}\textcolor{blue}{D_y}\mathcal{E}_{x};\Phi^v_x\right\rangle
$$
$$
+Bg\left\langle \textcolor{blue}{D}\textcolor{blue}{D_x}\mathcal{E}_{x}+2\textcolor{blue}{D}\textcolor{blue}{D_y}\mathcal{E}_{y};\Phi^v_y\right\rangle-\int_{\Gamma_n}Bg(\textcolor{blue}{D}\textcolor{blue}{D_x}+3\textcolor{blue}{D}\textcolor{blue}{D_y})\Phi^v\ds\frac{\partial \mathcal{E}}{\partial n}\partial \gamma.
$$

\begin{rem} We note that, during the simulation, when there is a steep gradient, we obtain a blow-ups. In order to avoid this problem, we need to change the bottom (make it smoother) or/and to get rid of the high order derivatives for the bottom as in \cite{Mit09}\label{Mit091}. That's why, we take into account in the sequel, the smoothness of the bottom and the fact that the derivatives of the bottom of order greater then one are neglected. 
\end{rem}
Thus, we will now deal with the following system:
\bg\left\{\label{derivBBMFullsimp}\begin{array}{rcl}
\Big\langle\partial_t\mathcal{E};\Phi^{\eta}\Big\rangle+b\left\langle \textcolor{blue}{D^2}\nabla\partial_t\mathcal{E};\nabla(\Phi^\eta)\right\rangle-\ds\int_{\Gamma_n}b\textcolor{blue}{D^2}\Phi^\eta\ds\frac{\partial (\partial_t\mathcal{E})}{\partial n}\partial \gamma &=&\mathcal{F}\left(\mathcal{E},\mathcal{U},\mathcal{V},\Phi^\eta\right)\\
 \Big\langle\partial_t\mathcal{U};\Phi^u\Big\rangle+d\left\langle\textcolor{blue}{D^2}\nabla\partial_t\mathcal{U};\nabla\Phi^u\right\rangle+d\left\langle2\textcolor{blue}{D}\nabla\textcolor{blue}{D}\cdot\nabla\partial_t\mathcal{U};\Phi^u\right\rangle-\ds\int_{\Gamma_n}d\textcolor{blue}{D^2}\Phi^u\ds\frac{\partial (\partial_t\mathcal{U})}{\partial n}\partial \gamma&=&\mathcal{G}\left(\mathcal{E},\mathcal{U},\mathcal{V},\Phi^u\right)\\
\Big\langle\partial_t\mathcal{V};\Phi^v\Big\rangle+d\left\langle\textcolor{blue}{D^2}\nabla\partial_t\mathcal{V};\nabla\Phi^v\right\rangle+d\left\langle2\textcolor{blue}{D}\nabla\textcolor{blue}{D}\cdot\nabla\partial_t\mathcal{V};\Phi^v\right\rangle-\ds\int_{\Gamma_n}d\textcolor{blue}{D^2}\Phi^v\ds\frac{\partial (\partial_t\mathcal{V})}{\partial n}\partial \gamma&=&\mathcal{H}\left(\mathcal{E},\mathcal{U},\mathcal{V},\Phi^v\right)
\end{array}\right.\ed
with
$$
\mathcal{F}\left(\mathcal{E},\mathcal{U},\mathcal{V},\Phi^\eta\right)=-\left\langle (\textcolor{blue}{D}+\textcolor{red}{\zeta}+\mathcal{E})\nabla\cdot(\mathcal{U};\mathcal{V})+(\textcolor{blue}{D_x}+\textcolor{red}{\zeta_x}+\mathcal{E}_x)\mathcal{U}+(\textcolor{blue}{D_y}+\textcolor{red}{\zeta_y}+\mathcal{E}_y)\mathcal{V}+\textcolor{red}{\zeta_t};\Phi^\eta\right\rangle
$$
$$
-\tilde{A}\left\langle 2\textcolor{blue}{D}\textcolor{blue}{D_x}\textcolor{red}{\zeta_{xt}}+2\textcolor{blue}{D}\textcolor{blue}{D_y}\textcolor{red}{\zeta_{yt}};\Phi^\eta\right\rangle-A\left\langle-2\textcolor{blue}{D}\textcolor{blue}{D_x}\textcolor{blue}{D_{y}}\mathcal{V}_{x}+2\textcolor{blue}{D}\textcolor{blue}{D_x^2}\mathcal{V}_{y};\Phi^\eta\right\rangle+A\left(\left\langle 2\textcolor{blue}{D^2}\textcolor{blue}{D_{x}}\mathcal{U}_{x}+\textcolor{blue}{D^2}\textcolor{blue}{D_{y}}\mathcal{V}_{x};\Phi^{\eta}_x\right\rangle+\left\langle \textcolor{blue}{D^2}\textcolor{blue}{D_{y}}\mathcal{U}_{x}\right.\right.
$$
$$
\left.\left.+\textcolor{blue}{D^2}\textcolor{blue}{D_{x}}\mathcal{U}_{y}+\textcolor{blue}{D^2}\textcolor{blue}{D_{x}}\mathcal{V}_{x}+2\textcolor{blue}{D^2}\textcolor{blue}{D_{y}}\mathcal{V}_{y};\Phi^{\eta}_y\right\rangle\right)-A\int_{\Gamma_n}\left((3\textcolor{blue}{D^2}\textcolor{blue}{D_{x}}+\textcolor{blue}{D^2}\textcolor{blue}{D_{y}})\Phi^\eta\ds\frac{\partial \mathcal{U}}{\partial n}+(\textcolor{blue}{D^2}\textcolor{blue}{D_{x}}+3\textcolor{blue}{D^2}\textcolor{blue}{D_{y}})\Phi^\eta\ds\frac{\partial \mathcal{V}}{\partial n}\right)\partial \gamma,
$$
$$
\mathcal{G}\left(\mathcal{E},\mathcal{U},\mathcal{V},\Phi^u\right)=-\left\langle g\left(I_d-B\left(2\textcolor{blue}{D_x^2}+\textcolor{blue}{D_y^2}\right)\right)\mathcal{E}_{x}+\mathcal{U}\mathcal{U}_{x}+\mathcal{V}\mathcal{V}_{x}-Bg\textcolor{blue}{D_x}\textcolor{blue}{D_{y}}\mathcal{E}_y-B\textcolor{blue}{D} \textcolor{red}{\zeta_{xtt}};\Phi^u\right\rangle+Bg\left\langle 2\textcolor{blue}{D}\textcolor{blue}{D_x}\mathcal{E}_{x};\Phi^u_x\right\rangle
$$
$$
+Bg\left\langle \textcolor{blue}{D}\textcolor{blue}{D_y}\mathcal{E}_{x}+\textcolor{blue}{D}\textcolor{blue}{D_x}\mathcal{E}_{y};\Phi^u_y\right\rangle-\int_{\Gamma_n}Bg(3\textcolor{blue}{D}\textcolor{blue}{D_x}+\textcolor{blue}{D}\textcolor{blue}{D_y})\Phi^u\ds\frac{\partial \mathcal{E}}{\partial n}\partial \gamma,
$$
and
$$
\mathcal{H}\left(\mathcal{E},\mathcal{U},\mathcal{V},\Phi^v\right)=-\left\langle -2Bg\textcolor{blue}{D_x}\textcolor{blue}{D_{y}}\mathcal{E}_x+\mathcal{U}\mathcal{U}_{y}+\mathcal{V}\mathcal{V}_{y}+g\left(I_d-2B\textcolor{blue}{D_y^2}\right)\mathcal{E}_{y}-B\textcolor{blue}{D} \textcolor{red}{\zeta_{ytt}};\Phi^v\right\rangle+Bg\left\langle \textcolor{blue}{D}\textcolor{blue}{D_y}\mathcal{E}_{x};\Phi^v_x\right\rangle
$$
$$
+Bg\left\langle \textcolor{blue}{D}\textcolor{blue}{D_x}\mathcal{E}_{x}+2\textcolor{blue}{D}\textcolor{blue}{D_y}\mathcal{E}_{y};\Phi^v_y\right\rangle-\int_{\Gamma_n}Bg(\textcolor{blue}{D}\textcolor{blue}{D_x}+3\textcolor{blue}{D}\textcolor{blue}{D_y})\Phi^v\ds\frac{\partial \mathcal{E}}{\partial n}\partial \gamma.
$$
\subsection{Time marching scheme}
Our method is based on an explicit second order Runge-Kutta scheme. For that, let us denote by $(\mathcal{E}^{n+1},\mathcal{U}^{n+1},\mathcal{V}^{n+1})$ and $(\mathcal{E}^n,\mathcal{U}^n,\mathcal{V}^n)$ the approximate values at time $t=t^{n+1}$ and $t=t^n$, respectively and by $\delta t$ the time step size. Then, owing to (\ref{derivBBMFullsimp}), the unknown fields at time $t=t^{n+1}$ are defined as the solution of the following system:

\bg\label{TBOUSW}\left\{\begin{array}{l}
\langle\mathcal{E}^{n+1};\Phi^\eta\rangle=\langle \mathcal{E}^n+\ds\frac{\mathcal{E}^{k1}+\mathcal{E}^{k2}}{2};\Phi^\eta\rangle,\\ \langle\mathcal{U}^{n+1};\Phi^u\rangle=\langle \mathcal{U}^n+\ds\frac{\mathcal{U}^{k1}+\mathcal{U}^{k2}}{2};\Phi^u\rangle,\\ \langle\mathcal{V}^{n+1};\Phi^v\rangle=\langle \mathcal{V}^n+\ds\frac{\mathcal{V}^{k1}+\mathcal{V}^{k2}}{2};\Phi^v\rangle,
\end{array}\right.\ed
where
\bg\label{PSIUVK1}\begin{array}{rcl}\left\langle\mathcal{E}^{k1};\Phi^\eta\right\rangle+b\left\langle  \textcolor{blue}{D^2}\nabla\mathcal{E}^{k1};\nabla(\Phi^\eta)\right\rangle-\ds\int_{\Gamma_n}b \textcolor{blue}{D^2}\Phi^\eta\ds\frac{\partial (\mathcal{E}^{k1})}{\partial n}\partial \gamma&=&\delta t\cdot\mathcal{F}\left(\mathcal{E}^n,\mathcal{U}^n,\mathcal{V}^n,\Phi^\eta\right),\\ 
\left\langle\mathcal{U}^{k1}+2d \textcolor{blue}{D}\nabla \textcolor{blue}{D}\cdot\nabla\mathcal{U}^{k1};\Phi^u\right\rangle+d\left\langle \textcolor{blue}{D^2}\nabla\mathcal{U}^{k1};\nabla\Phi^u\right\rangle-\ds\int_{\Gamma_n}d \textcolor{blue}{D^2}\Phi^u\ds\frac{\partial (\mathcal{U}^{k1})}{\partial n}\partial \gamma&=&\delta t\cdot\mathcal{G}\left(\mathcal{E}^n,\mathcal{U}^n,\mathcal{V}^n,\Phi^u\right),\\
\left\langle\mathcal{V}^{k1}+2d \textcolor{blue}{D}\nabla \textcolor{blue}{D}\cdot\nabla\mathcal{V}^{k1};\Phi^v\right\rangle+d\left\langle \textcolor{blue}{D^2}\nabla\mathcal{V}^{k1};\nabla\Phi^v\right\rangle-\ds\int_{\Gamma_n}d \textcolor{blue}{D^2}\Phi^v\ds\frac{\partial (\mathcal{V}^{k1})}{\partial n}\partial \gamma&=&\delta t\cdot\mathcal{H}\left(\mathcal{E}^n,\mathcal{U}^n,\mathcal{V}^n,\Phi^v\right)
\end{array}\ed
and
\bg\label{PSIUVK2}\hspace{-1cm}\begin{array}{rcl} \left\langle\mathcal{E}^{k2};\Phi^\eta\right\rangle+b\left\langle  \textcolor{blue}{D^2}\nabla\mathcal{E}^{k2};\nabla(\Phi^\eta)\right\rangle-\ds\int_{\Gamma_n}b \textcolor{blue}{D^2}\Phi^\eta\ds\frac{\partial (\mathcal{E}^{k2})}{\partial n}\partial \gamma&=&\delta t\cdot\mathcal{F}\left(\mathcal{E}^n+\mathcal{E}^{k1},\mathcal{U}^n+\mathcal{U}^{k1},\mathcal{V}^n+\mathcal{V}^{k1},\Phi^\eta\right),\\
\left\langle\mathcal{U}^{k2}+2d \textcolor{blue}{D}\nabla \textcolor{blue}{D}\cdot\nabla\mathcal{U}^{k2};\Phi^u\right\rangle+d\left\langle \textcolor{blue}{D^2}\nabla\mathcal{U}^{k2};\nabla\Phi^u\right\rangle-\ds\int_{\Gamma_n}d \textcolor{blue}{D^2}\Phi^u\ds\frac{\partial (\mathcal{U}^{k2})}{\partial n}\partial \gamma&=&\delta t\cdot\mathcal{G}\left(\mathcal{E}^n+\mathcal{E}^{k1},\mathcal{U}^n+\mathcal{U}^{k1},\mathcal{V}^n+\mathcal{V}^{k1},\Phi^u\right),\\
\left\langle\mathcal{V}^{k2}+2d \textcolor{blue}{D}\nabla \textcolor{blue}{D}\cdot\nabla\mathcal{V}^{k2};\Phi^v\right\rangle+d\left\langle \textcolor{blue}{D^2}\nabla\mathcal{V}^{k2};\nabla\Phi^v\right\rangle-\ds\int_{\Gamma_n}d D^2\Phi^v\ds\frac{\partial (\mathcal{V}^{k2})}{\partial n}\partial \gamma&=&\delta t\cdot\mathcal{H}\left(\mathcal{E}^n+\mathcal{E}^{k1},\mathcal{U}^n+\mathcal{U}^{k1},\mathcal{V}^n+\mathcal{V}^{k1},\Phi^v\right).
\end{array}\ed

\section{Mesh generation and initial data}\label{meshgeninitdata}
In this section, we present the technique used in order to generate a mesh using a photo of the Mediterranean sea then using an imported xyz bathymetry. Also, we will present a way in order to obtain the initial data and explain the special adapt mesh technique that we will use in our numerical simulation.

\subsection{Mesh generated using a photo}\label{meshphoto}
We present here the method to build a mesh from a photo inspired from a \freefem script made by F. Hecht \cite{Hecht11}\label{Hecht110} and another one made by O. Pantz \cite{Pan11}\label{Pan110}.

Owing to Google Earth$\circledR$ and for better resolution, we take severals parts of the Mediterranean sea (cf. Figure \ref{Mediterannee_coupe}) that are subsequently assembled using Photoshop$\circledR$ to obtain a complete picture of the Mediterranean sea (cf. Figure \ref{Mediterannee_ALL}).
\begin{figure}[!h]
\begin{center}
    \includegraphics[height=4cm,width=6cm]{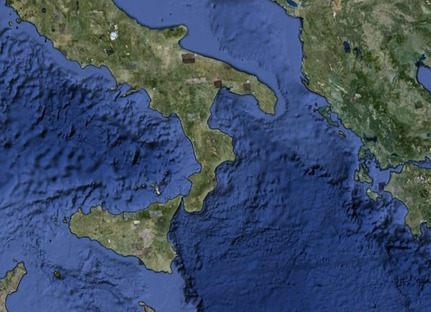} $\qquad$\includegraphics[height=4cm,width=6cm]{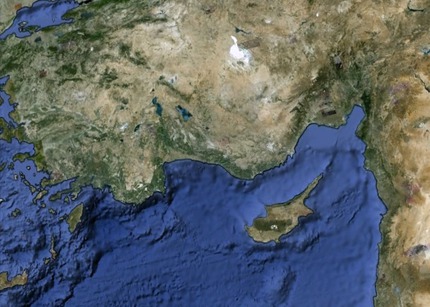}
\caption{The pictures of two parts of the Mediterranean sea. \label{Mediterannee_coupe} }
\end{center}
\end{figure}
\clearpage
\begin{figure}[!h]
\begin{center}
 \includegraphics[height=7cm,width=14.5cm]{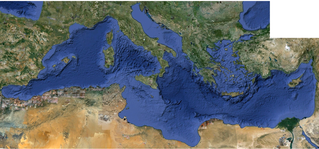}
\caption{The Mediterranean sea after assembly with Photoshop$\circledR$. \label{Mediterannee_ALL} }
\end{center}
\end{figure}
Using Photoshop$\circledR$, we can also eliminate the land areas that circumvent the Mediterranean sea (cf. Figure \ref{Mediterannee_ALL_wet}). We note that we must smoothen the borders in Photoshop$\circledR$ and put the black color inside our domain and the white color outside.\\
\begin{figure}[!h]
\begin{center}
    \includegraphics[height=7cm,width=14.5cm]{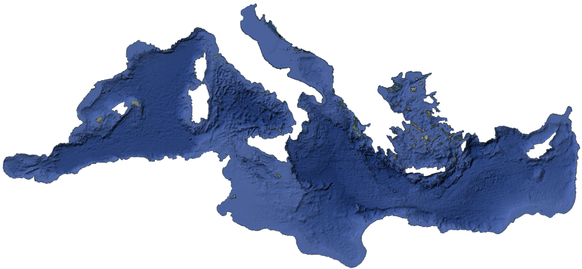}
    \end{center}
\caption{The Mediterranean sea. \label{Mediterannee_ALL_wet} }
\end{figure}

Then we convert the jpg photo to a pgm photo which can be read by \freefem using in a terminal window :
\begin{lstlisting}[firstnumber=last]
convert Medit_sea.jpg Medit_sea.pgm
\end{lstlisting}
In order to generate the mesh of the Mediterranean sea domain, we read the pgm file using
\begin{lstlisting}[firstnumber=last]
load "ppm2rnm"
load "isoline"
real[int,int] f1("Medit_sea.pgm");
int nx = f1.n, ny=f1.m;
mesh Sh=square(nx-1,ny-1,[(nx-1)*(x),(ny-1)*(1-y)]);
fespace Vh(Sh,P1);	Vh fxy;
fxy[]=f1;
fxy=(fxy>=0.5)-(fxy<0.5);// to get value between -1 and 1
real[int,int] Curves(3,1);
int[int] be(1);
int nc=isoline(Sh,fxy,iso=0.,close=0,Curves,beginend=be,smoothing=.005,ratio=0.1);
\end{lstlisting}
The function obtained from the pgm file has values between $0$ and $1$, where the value of $0.5$ represents the contour between two different colors. We note that, we can regularize this contour (thanks to O. Pantz \cite{Pan11}\label{Pan111}), before using the isoline function which computes the number of closed curves of our image, by solving :
\bg\label{reg}\ds\int_\Omega\left(2\varepsilon\nabla(\delta u)\nabla(\phi)+\ds\frac{4}{\varepsilon}u^2\delta u\phi+\alpha \delta u\phi\right)+\ds\int_\Omega\left(\varepsilon\nabla(u)\nabla(\phi)-\ds\frac{2}{\varepsilon}(1-u^2)u\phi\right)=\ds\int_\Omega\alpha(fxy-u)\phi,\ed
where in our example, we take $\varepsilon=1,\ \alpha=1$ and $\delta x=1.5$. We also note that when $\varepsilon$ is close to zero, the solution $u$ takes the values equal to $1$ and $-1$ and $\alpha$ sets the balance between length of the curve and fitting the actual interface: As $\alpha$ increase, the approximation become better.\\
After regularizing, we update $u=0$ till $u=u+\delta u$, and we adapt the mesh around the contour (cf. Figure \ref{Mediterannee_ALL_mesh}), by using 
\begin{lstlisting}[firstnumber=last]
Sh=adaptmesh(Sh,u,nbvx=1e7);
\end{lstlisting}
\begin{figure}[!h]
\begin{center}
\includegraphics[height=7cm,width=14.5cm]{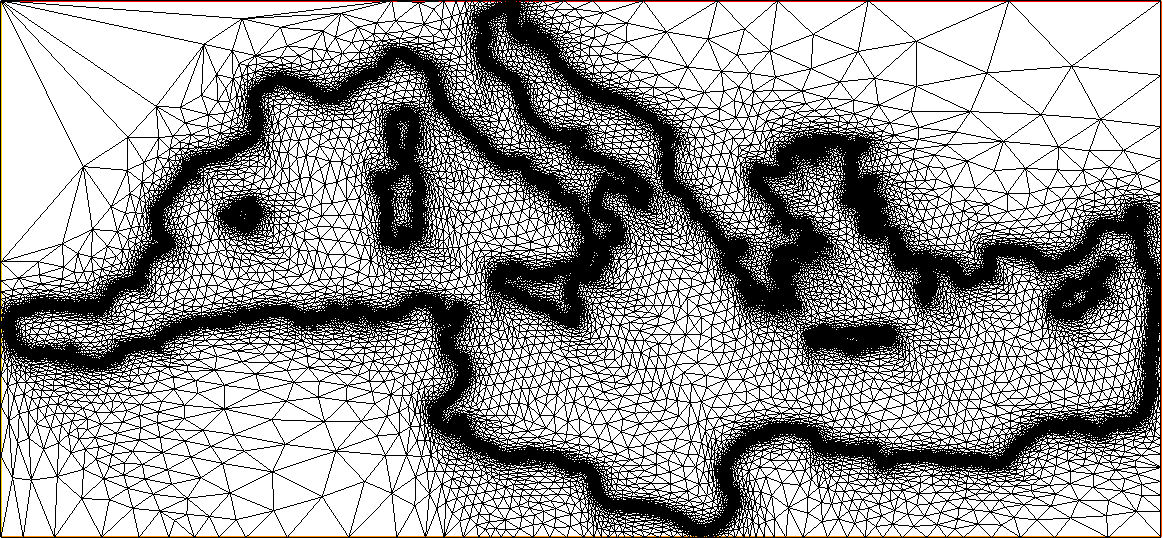}
\end{center}
\caption{The adapted mesh around the contour of the Mediterranean sea. \label{Mediterannee_ALL_mesh} }
\end{figure}
then, we solve again the regularizing problem, then we update the function $u=u+\delta u$ and finally we interpolate the solution on the initial mesh (cf. Figure \ref{Mediterannee_ALL_interpolate}) before using isoline :
\begin{lstlisting}[firstnumber=last]
mesh Shinit=Sh;	fespace Vhinit(Shinit,P1);
Vhinit fxyinit;
fxyinit=u;
int nc=isoline(Shinit,fxyinit,iso=0.,close=0,Curves,beginend=be,smoothing=.005,ratio=0.1);
\end{lstlisting}

Now, we show in Figure \ref{Crete_mesh}, the mesh created around Crete island, where we see in the top left, the mesh after first regularization and in the top right, the function $fxy$ after interpolation, in the down left, three isoline level $-1$, $0$ and $1$ of the $fxy$ function, where here the mesh generated at the $0$ level is shown in the down right. The complete script is written in \cite{Sad12e}\label{Sad12e0}.\\
We note that, we take into account the area of the Mediterranean sea, which is almost $2.5$ million $km^2$ and we did a scale for our final mesh :
\begin{lstlisting}[firstnumber=last]
real Areasea=2.5e6;//area of the Mediterranean sea in Km^2
real scale=sqrt(Areasea/Th.area);
Th=movemesh(Th,[x*scale,y*scale]);
\end{lstlisting}
\clearpage
\begin{figure}[!h]
\begin{center}
\includegraphics[height=7cm,width=14.5cm]{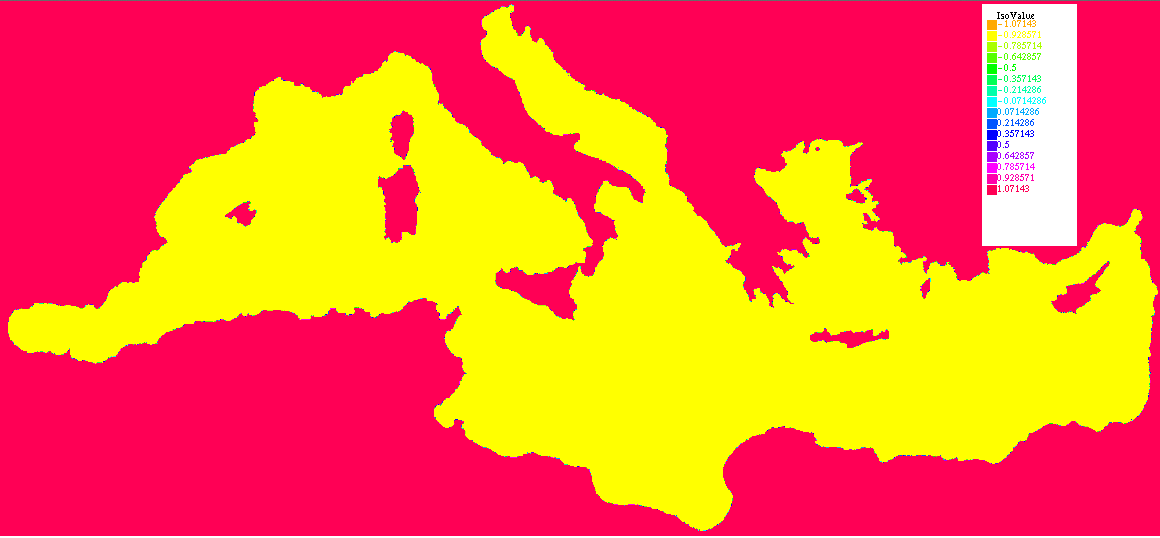}
\end{center}
\caption{The function $fxy$ of the contour of Mediterranean sea after interpolation.\label{Mediterannee_ALL_interpolate} }
\end{figure}
\begin{figure}[!h]
\begin{center}
    \includegraphics[height=5cm,width=8cm]{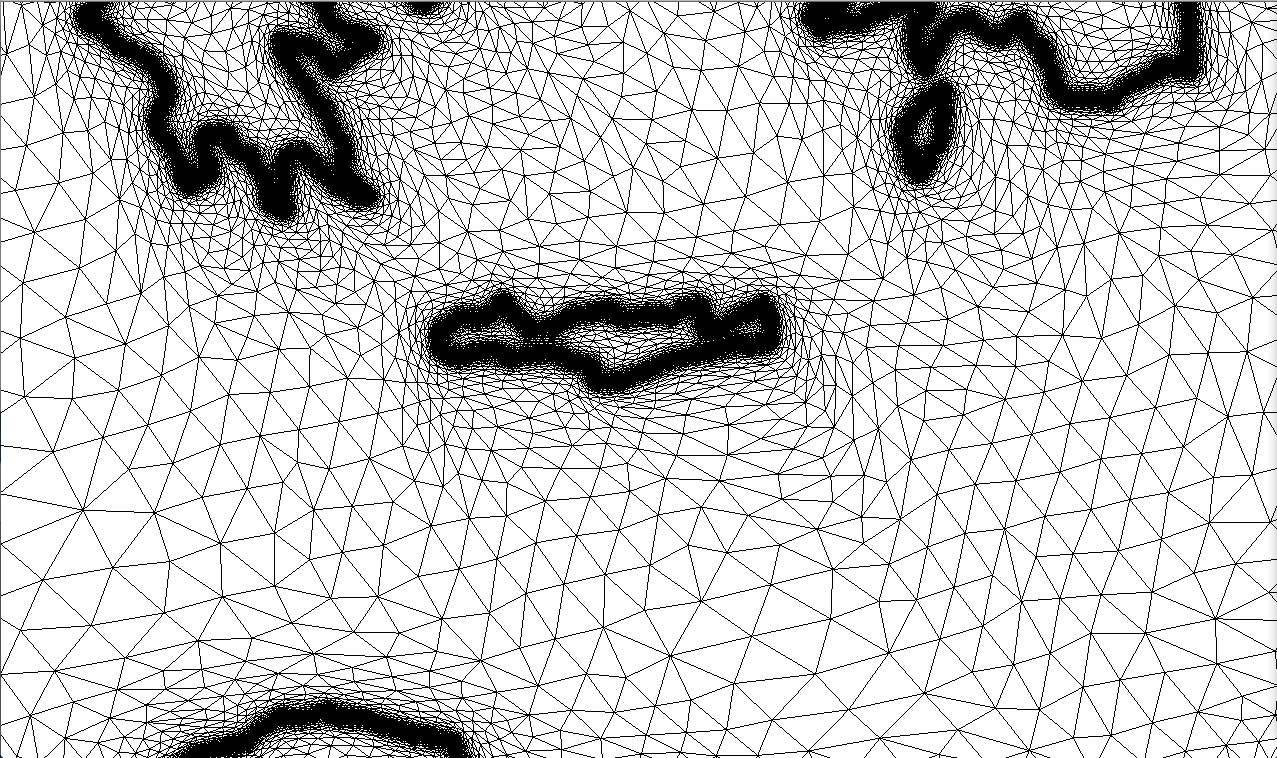} $\qquad$\includegraphics[height=5cm,width=8cm]{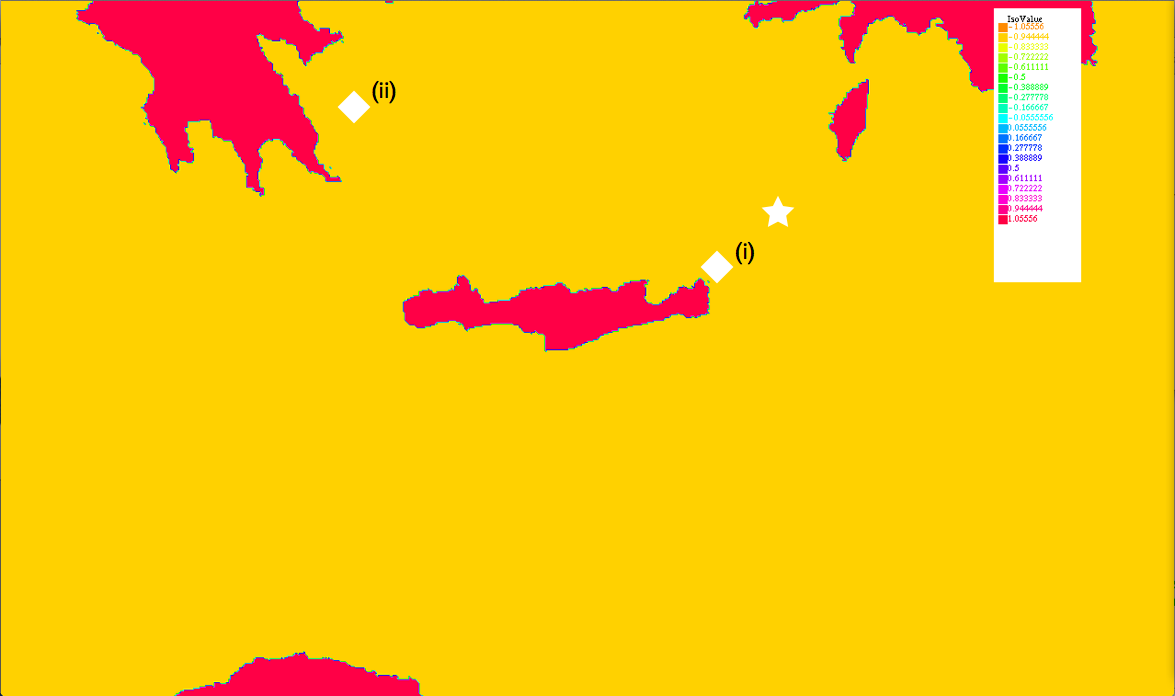}
        \includegraphics[height=5cm,width=8cm]{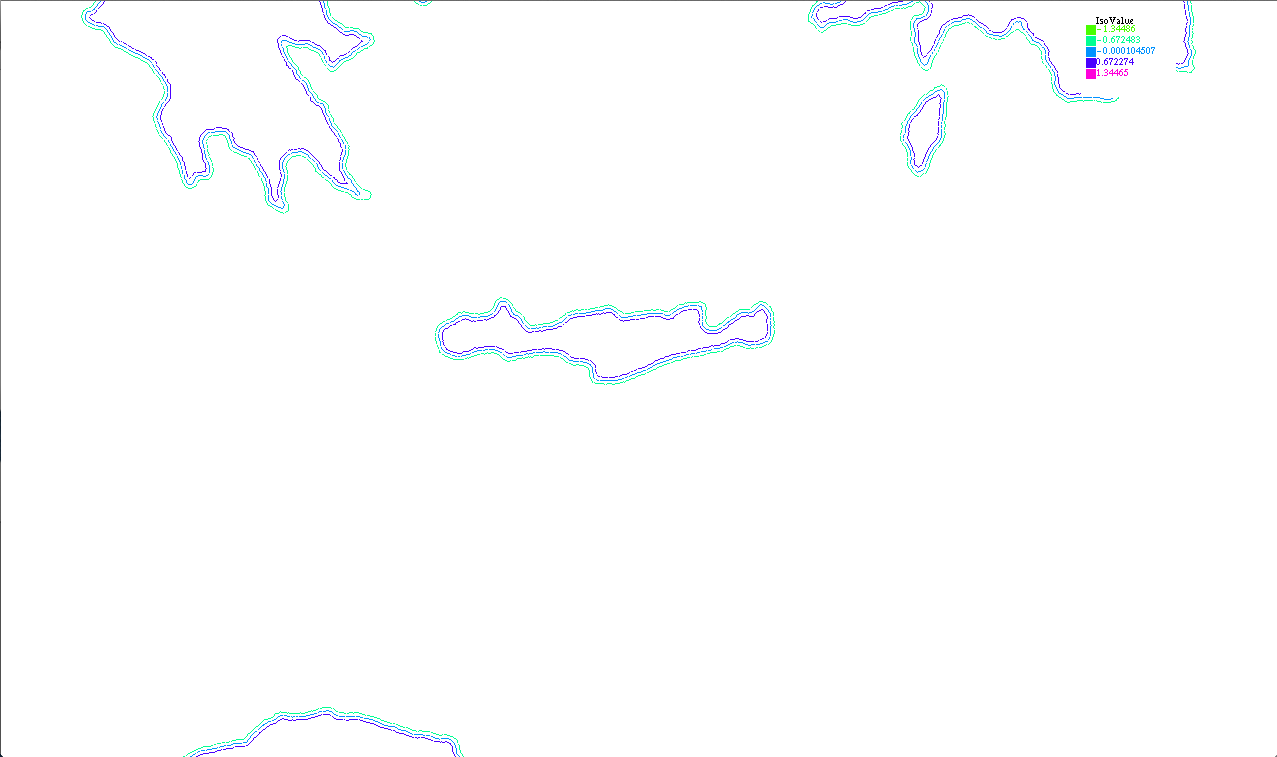} $\qquad$\includegraphics[height=5cm,width=8cm]{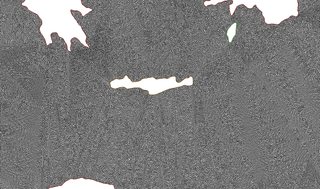}
\caption{Creation of the mesh around Crete island in the Mediterranean sea (top left: mesh after first regularization, top right: solution after interpolation, $\diamond$ : wave gauge, $\star$ : epicenter, down left: isoline of the contour, down right: mesh generated).\label{Crete_mesh} }
\end{center}
\end{figure}

\subsection{Mesh generated using an imported xyz bathymetry}\label{secgenmeshxyz}
In order to consider more realistic case, this means that, we now take into account the bathymetry near Java island which can be downloaded from the NOAA\footnote{\url{https://maps.ngdc.noaa.gov/viewers/wcs-client/}} website, (cf. Figure \ref{batdat}), we also use \freefem to generate the mesh of the area where the amplitude is zero.\\
\begin{figure}[h]
\begin{center}
    \includegraphics[height=6cm,width=14cm]{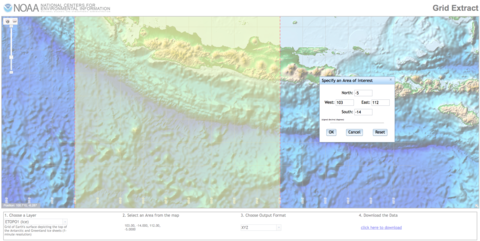}
\end{center}
\caption{\label{batdat}Importation of bathymetry datas through the NOAA website.}
\end{figure}

We can read the xyz file (cf. Figure \ref{bottom_Java}, left), using this script :
\begin{lstlisting}[firstnumber=last]
mesh Sh=triangulate("bathymetry_Java");
fespace Vh(Sh,P1);	Vh fxy;
{ ifstream file(filename);real xx,yy;
	for(int i=0;i<fxy.n;i++)
		file >> xx >>yy >> fxy[][i];
}
\end{lstlisting}
We can smoothen the bathymetric data by solving :

\bg\label{smbot}\ds\int_\Omega\left(\beta u\phi+\nabla(u)\nabla(\phi)\right)+\ds\int_\Gamma\ds\frac{\partial u}{\partial n}\partial \gamma=\ds\int_\Omega\beta fxy\phi,\qquad \ds\frac{\partial u}{\partial n}=0.\ed

In our code, we take $\beta=5.e3$, in order to build the mesh to get rid of smallest islands, then, in order to obtain the mesh only around the sea, which is limited by the zero level of the amplitude, we use :
\begin{lstlisting}[firstnumber=last]
fxy=(fxy>=0.5)-(fxy<0.5);// to get value between -1 and 1
\end{lstlisting}
and then proceeding similarly, as above for the mesh generation using a photo, in order to obtain the mesh of our domain, by taking $\varepsilon=1.e-2$ and $\alpha=1.e2$ for the regularization phenomena (\ref{reg}).

\begin{rem}
We note that our method takes into account the different label for each part of the boundary, which facilitates the use of different types of boundary condition.
\end{rem}

\begin{rem}
For all simulation with bathymetry, we use $\beta=2.e1$ in (\ref{smbot}) to smoothen the initial bathymetry after generation of the mesh (cf. Figure \ref{bottom_Java}, right) in order to ensure the stability of the numerical method, we also note that in order to be in a big deep water wave regime for sBBM system we change the depth close to the shoreline to 100 m.
\end{rem}

\begin{figure}[!h]
\begin{center}
\begin{tabular}{m{7cm}m{7cm}}
    \includegraphics[height=5cm,width=7cm]{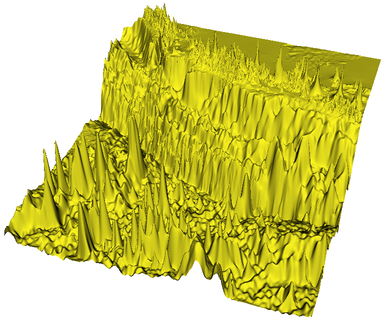} &\includegraphics[height=5cm,width=7cm]{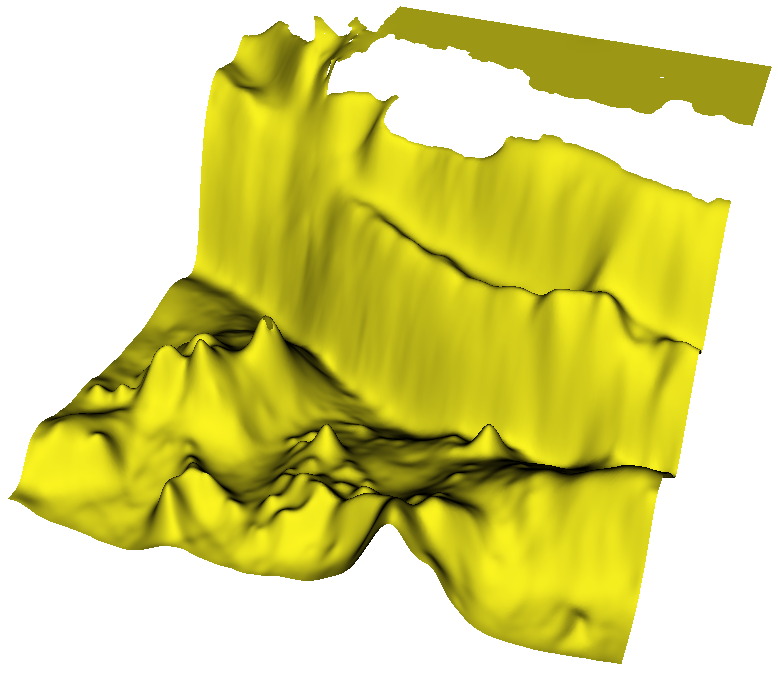}
\end{tabular}    
\caption{Left : Bathymetry downloaded from the NOAA website, (min $=-7239m$ and max $=3002m$).\\
Right : smoothed bathymetry with $\beta=2.e1$ in (\ref{smbot}), (min $=-6207m$ and max $=-100m$).\label{bottom_Java} }
\end{center}
\end{figure}
\begin{rem}
The bathymetry data downloaded from the NOAA website are in degree coordinate and we need to convert them to meter. So, on the first hand, we must know the degree of Latitude (south and north) and of Longitude (west and east) of our domain where we can deduce the Latitude $lat0=.5(lat_{south}+lat_{north})$ and the Longitude $long0=.5(long_{west}+long_{east})$. On the other hand, we must take into account the spherical shape of the Earth, even if it does not play significant role because of the small spatial scale of the experiments. So, we know that the radius of the Earth near the equator is $R_{equator}=6378,137$ Km, and near to the pole $R_{pole}=6356,752$ Km, thus the radius of our domain equals to:
$$
R=\sqrt{\ds\frac{\big(R_{equator}^2\cos(lat0\cdot \pi/180)\big)^2+\big(R_{pole}^2\sin(lat0\cdot \pi/180)\big)^2}{\big(R_{equator}\cos(lat0\cdot \pi/180)\big)^2+\big(R_{pole}\sin(lat0\cdot \pi/180)\big)^2}}.
$$
So, we move the mesh of our domain using the following translation ($\mbox{coefl0}=\pi R/180$):
$$
[x;y]\longrightarrow [(x-lon0)\cos(\pi y/180)\mbox{coefl0};(y-lat0)\mbox{coefl0}].
$$
Therefore, we need to move the bottom from the original reference downloaded to the new mesh, which is easy to do in \freefem by writing this script :
\begin{lstlisting}[firstnumber=last]
Vh D=fxy;
real coefl0=pi*R/180.;
Sh=movemesh(Sh,[(x-lon0)*cos(y*pi/180.)*coefl0,(y-lat0)*coefl0]);
Vh tmp;	// define a temporary function
tmp=D[];	// save the value
D=0;		// to change the FEspace and mesh associated with u
D[]=tmp;	// set the value of u without any mesh update
\end{lstlisting}
\end{rem}

\subsection{Mesh adaptation technique}
We introduce here a special mesh adaptation technique since some computation domains are huge as in the case of the Mediterranean sea with $1.7e6$ triangles ($8.7e5$ degree of freedom) and our initial solution is concentrated in a small domain, a circle $\mathcal{C}(O,R)$ or a rectangle $[aa,bb]\times[cc,dd]$. So we build the mesh {\ttfamily{Th}} of the small rectangle or the circle by doing a trunk to the initial full mesh {\ttfamily{Thinit}} respecting that the function equals to $1$ inside the rectangle and $0$ outside as in : 
 \begin{lstlisting}[firstnumber=last]
Th=trunc(Thinit,(1.*(x <= bb & x >= aa) *(y <= dd & y >= cc))>.5,label=0);
\end{lstlisting}
or
 \begin{lstlisting}[firstnumber=last]
Th=trunc(Thinit,(1.*((x -xO)^2+(y-yO)^2<=R^2))>.5,label=0);
\end{lstlisting}
then we compute the initial solution {\ttfamily{uold=uadapt}} in this domain. Using the keyword \textcolor{red}{\ttfamily{boundingbox}} in \freefem, we obtain the limit min max of {\ttfamily{Thold=Th}} on $x$ and $y$ direction, in which we add {\ttfamily{epsadapt}} from each side in order to build the new rectangle {\ttfamily{Thnew}} that contains {\ttfamily{Thold}}, then using the keyword \textcolor{red}{\ttfamily{interpolate}} in \freefem, the old FEspace in $\mathbb{P}_1$ {\ttfamily{Vhold}} and the new FEspace in $\mathbb{P}_1$ {\ttfamily{Vhnew}}, we interpolate {\ttfamily{uold}} to {\ttfamily{unew}} in the new domain. We smooth the function obtained from {\ttfamily{abs(unew)>=erradapt}} using :
\bg\label{smbotdir}\ds\int_\Omega\left(\mbox{{\ttfamily{smoothadapt}} } \cdot \mbox{ {\ttfamily{usmadapt}}}\cdot\phi+\nabla(\mbox{{\ttfamily{usmadapt}}})\nabla(\phi)\right)=\ds\int_\Omega{\mbox{\ttfamily{smoothadapt}} }\cdot\left( |\mbox{{\ttfamily{unew}}}|\geq \mbox{{\ttfamily{erradapt}}}\right)\cdot\phi,\ed
with zero Dirichlet boundary condition on the boundary label 0 of {\ttfamily{Thnew}} in order to obtain {\ttfamily{usmadapt}}. And, at the end, we trunk {\ttfamily{Thnew}} respecting the function {\ttfamily{usmadapt>=isoadapt}}. We can put all the previous detail for the special mesh adaption in a \textcolor{red}{\ttfamily{macro}} such as :
 \begin{lstlisting}[firstnumber=last]
macro AdaptGS(Th,uadapt,erradapt,isoadapt,smoothadapt,epsadapt,Thinit,output,wv)
	mesh Thold=Th;
	fespace Vhold(Thold,P1);
	Vhold uold=uadapt;
	real[int] bbmM(4);
	boundingbox(Thold,bbmM);
	mesh Thnew=trunc(Thinit,(1.*(x <= (bbmM[1]+epsadapt) & x >= (bbmM[0]-epsadapt)) *(y <= (bbmM[3]+epsadapt) & y >= (bbmM[2]-epsadapt)))>erradapt,label=0);
	fespace Vhnew(Thnew,P1);
	Vhnew usmadapt,vsmadapt,unew;
	matrix BSinterp=interpolate(Vhnew,Vhold,inside=true);
	unew[]=BSinterp*uold[];
	solve smoother(usmadapt,vsmadapt)=int2d(Thnew)(smoothadapt*usmadapt*vsmadapt+grad(usmadapt)'*grad(vsmadapt))-int2d(Thnew)(smoothadapt*(abs(unew)>=erradapt)*vsmadapt)+on(0,usmadapt=0.);
	if(output) plot(usmadapt,fill=true,dim=2,value=true,wait=wv,cmm="usmadapt");
	if(output) {Vhnew unewp=(usmadapt>=isoadapt); plot(unewp,fill=true,dim=2,value=true,wait=wv,cmm="zone to cut");}
	Th=trunc(Thnew,usmadapt>=isoadapt,label=0);
//
\end{lstlisting}
\begin{rem} We trunk always from the initial full mesh, in this case, we keep the original vertices of the mesh throughout the simulation, and also we keep the original label of the boundary and we put 0 for the label of the rest of the boundary domain. We also note that we need to interpolate all the variables of any kind of FEspace from the old mesh to the newest one using the keyword \textcolor{red}{\ttfamily{interpolate}} in \freefem. We remark also that we must choose the parameters for the {\ttfamily{AdaptGS}} in order to obtain the value of {\ttfamily{erradapt}} on the boundary of {\ttfamily{Th}} for the function {\ttfamily{uadapt}}. In addition, we use a reflective boundary condition (BC) on label 0, {\it i.e.} zero Neumann BC for $\eta$ and zero Dirichlet BC for $\mathcal{V}$, cause our sBBM system gives artificial numerical explosion on the boundary if we do not use any BC or if we use only Neumann BC for $\eta$ and $\mathcal{V}$.
\end{rem}

\subsection{Initial data}
Inspiring from \cite{DMGD13,Mit09}\label{DMGD130}\label{Mit092}, {\it Tsunami} waves are generated by a deformation of the bottom due to an Earthquake, which may be approximated by Okada's formula's \cite{Oka85,Oka92}\label{Oka851}\label{Oka921} in the case of the so called dip-slip dislocation, where the vertical component of displacement vector $\mathcal{O}(x,y)$, is given by the following formulas in Chinnery's notation, cf \cite{DD07,Oka85}\label{DD071}\label{Oka852}
$$
f(\xi,\eta)\ds\left|\left|=f\left(\xi,p\right)-f\left(\xi,p-W\right)-f\left(\xi-L,p\right)+f\left(\xi-L,p-W\right)\right.\right.,
$$ 
$$
\left.\left.\mathcal{O}(x,y)=-\frac{U}{2\pi}\left(\frac{\tilde{d}q}{R(R+\xi)}+\sin\delta\arctan\frac{\xi\eta}{qR}-I\sin\delta\cos\delta\right)\right|\right|.
$$
where 
$$
\xi=(x-x_0)\cos\phi+(y-y_0)\sin\phi,\ Y=-(x-x_0)\sin\phi+(y-y_0)\cos\phi,
$$
$$
p=Y\cos\delta+d\sin\delta,\ q=Y\sin\delta-d\cos\delta,
$$
$$
\tilde{y}=\eta\cos\delta+q\sin\delta,\ \tilde{d}=\eta\sin\delta-q\cos\delta,
$$
$$
R^2=\xi^2+\eta^2+q^2=\xi^2+\tilde{y}^2+\tilde{d}^2,\ X^2=\xi^2+q^2
$$
and 
$$
I=\left\{
\begin{array}{ll}
\ds\frac{\mu}{\lambda+\mu}\ds\frac{2}{\cos\delta}\arctan\ds\frac{\eta(X+q\cos\delta)+X(R+X)\sin\delta}{\xi(R+X)\cos\delta}&\mbox{ if }\cos\delta\neq 0,\\  
\ds\frac{\mu}{\lambda+\mu}\ds\frac{\xi\sin\delta}{R+\tilde{d}}&\mbox{ if }\cos\delta= 0.
\end{array}\right.
$$
Here, $W$ and $L$ are the width and the length of the rectangular fault, $(x,y)$ are the points where we computes displacements, $(x_0,y_0)$ is the epicenter,\ $d=\mbox{fault depth}(x_0,y_0)+W\sin\delta$, $\delta$ is the dip angle,\ $\theta$ is the rake angle,\ $D$ is the Burger's vector,\ $U=|D|\sin\theta$ is the slip on the fault,\ $\phi$ is the strike angle which is measured conventionally in the counter-clockwise direction from the North (cf. Figure \ref{Okada_domain_sol} (left)),\ $\mu,\ \lambda$ are the Lam\'e constants derived from elastic-wave velocities : $\lambda=\rho_c\left(V_P^2-V_S^2\right)$ and $\mu=\rho_cV_S^2$, where $\rho_c$ is the crust density, $V_P$ is the compressional-wave (P-wave) velocity, $V_S$ is the shear-wave (S-wave) velocity.
\begin{rem}
We can download the script which computes co-seismic displacements according to the classical Okada solution $\mathcal{O}(x,y)$ from the following link \url{http://www.denys-dutykh.com/downloads.php}.
\end{rem}
We will distinguish here the two cases for the mechanisms of the dynamics of {\it Tsunami} wave generation as in \cite{Mit09}\label{Mit093} : the {\it passive} generation and the {\it active} generation.
\subsubsection*{{\it Passive} generation :}
In order to compute the initial data for $\eta(x,y,0)=\mathcal{O}(x,y)$ in meters (cf. Figure \ref{Okada_domain_sol} (right)), $V(x,y,0)=0$ which is referred to as a {\it passive} generation of a {\it Tsunami} wave near Java island, using our mesh adaptive technique, we will use the fact that the solution is concentrated in the small rectangle $[x_0-3.2W;x_0+1.2W]\times[y_0-L;y_0+L]$ where $L=100$ Km, $W=50$ Km, $\delta=10.35\degree$, $\phi=288.94\degree$, $\theta=95\degree$, $U=2$ m, $\rho_c=2700\mbox{ Kg/m}^3$, $V_P=6000$ m/s, $V_S=3400$ m/s, $(x_0;y_0)=(107.345\degree,-9.295\degree)$ and the fault depth$(x_0;y_0)=10$ Km. All these geophysical parameters can be downloaded from this website \href{https://Earthquake.usgs.gov/archive/product/finite-fault/usp000ensm/us/1486510367579/web/p000ensm.param}{https://Earthquake.usgs.gov}.

Therefore, we build the mesh of the small rectangle by doing a \textcolor{red}{\ttfamily{trunc}} to the initial full mesh respecting that the function equals to $1$ inside the rectangle and $0$ outside as in :
 \begin{lstlisting}[firstnumber=last]
real aa=x0-3.2*W,bb=x0+1.2*W,cc=y0-L,dd=y0+L;
Th=trunc(Thinit,(1.*(x <= bb & x >= aa) *(y <= dd & y >= cc))>.5,label=0);
\end{lstlisting}
\begin{figure}[!h]
\begin{center}
\begin{tabular}{m{7cm}m{7cm}}
    \includegraphics[height=3cm,width=7cm]{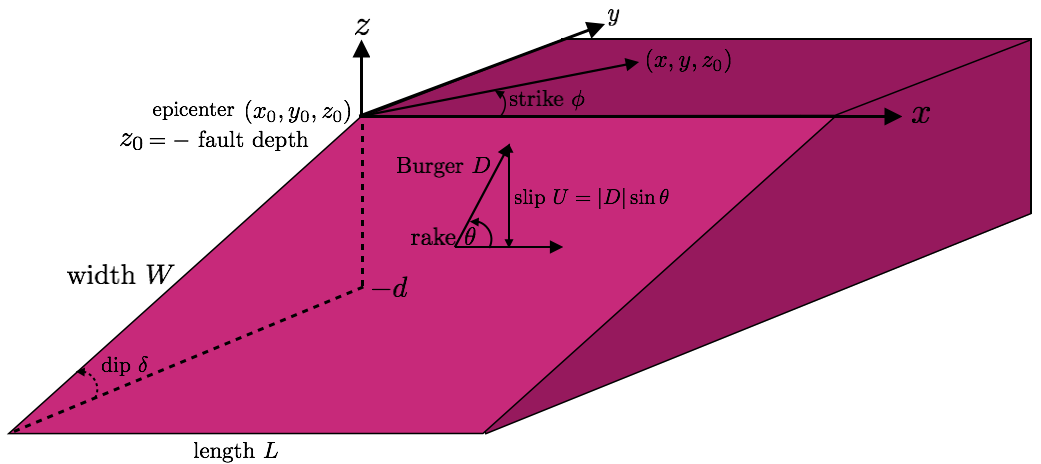} &\includegraphics[height=5cm,width=7cm]{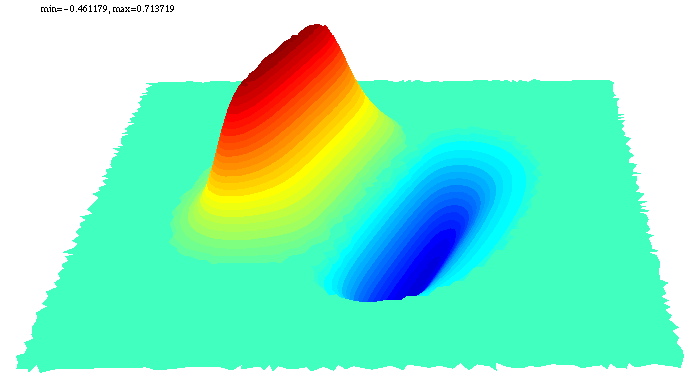}
\end{tabular}    
\caption{Geometry of the source model (left) and the initial solution for $\eta$ (right, min$=-0.46$ m, max$=0.73$ m ).\label{Okada_domain_sol} }
\end{center}
\end{figure}

\subsubsection*{{\it Active} generation :}
For a more realistic case as in the Java 2006 event, we use the {\it active} generation in order to model the generation of a {\it Tsunami} wave as in \cite{DMGD13,DMCS12}\label{DMGD131}\label{DMCS120}. In this case we consider zero initial conditions for both the free surface elevation and the velocity field, and assume that the bottom is moving in time. This case may be described by considering the bottom motion formula :
$-h(x,y,t)=-\textcolor{blue}{D}(x,y)-\textcolor{red}{\zeta}(x,y,t)$ with 
$$\textcolor{red}{\zeta}(x,y,t)=\ds\sum_{i=1}^{Nx\cdot Ny}\mathcal{H}(t-t_i)\cdot \left(1-e^{-\alpha(t-t_i)}\right)\cdot\mathcal{O}_i(x,y),$$
where $Nx$ sub-faults along strike and $Ny$ sub-faults down the dip angle, $\mathcal{H}(t)$ is the Heaviside step function and $\alpha=\frac{\log(3)}{t_r}$, where $t_r=8$ s is the rise time. We choose here an exponential scenario, but in practice, various scenarios could be used (instantaneous, linear, trigonometric, etc) and could be found in \cite{DMGD13,DMCS12,Ham72}\label{Ham720}\label{DMGD132}\label{DMCS121}. 
\begin{rem}
Parameters such as sub-fault location $(x_i,y_i)$, depth $d_i$, slip $U$ and rake angle $\theta$ for each segment, given in Table 3 of the paper \cite{DMGD13}\label{DMGD133}, can be downloaded from this website \href{https://Earthquake.usgs.gov/archive/product/finite-fault/usp000ensm/us/1486510367579/web/p000ensm.param}{https://Earthquake.usgs.gov}.\\
In this file, we remark that the fault's plane is conventionally divided into $Nx = 21$ sub-faults along strike and $Ny = 7$ sub-faults down the dip angle, leading to a total number of $Nx \times Ny = 147$ equal segments.
\end{rem}
For our special adapt mesh technique, since the fault plane is considered to be the rectangle with vertices located at ($109.20508\degree$ (Lon),$-10.37387\degree$ (Lat)), ($106.50434\degree$ (Lon),$-9.45925\degree$ (Lat)), ($106.72382\degree$ (Lon),$-8.82807\degree$ (Lat)) and ($109.42455\degree$ (Lon),$-9.74269\degree$ (Lat)), we will consider that our bottom displacement is concentrated on the big rectangle which is equidistant of $1\degree$ from each side of the initial fault plane as in Figure \ref{Fault_Okada_domain_mesh_sol} (left), then we compute each Okada solution $\mathcal{O}_i$ on a circle of center $(x_i-10m,y_i-10m)$ and of radius $6\max(L,W)$ and at then end all the Okada solution will be interpolated on the big rectangle before starting to compute the vertical displacement of the bottom $\textcolor{red}{\zeta}(x,y,t)$, in Figure \ref{Fault_Okada_domain_mesh_sol} (right) we plot $\mathcal{O}_{14}$. For the computation of $\textcolor{red}{\zeta}(x,y,t)$, we start the mesh by a circle of center $(x_c-5m,y_c-5m)$ and of radius $4\max(L,W)$ and we adapt the mesh each 3 iterations {\it i.e.} each $6$ s by using the following value for the adapt mesh {\ttfamily{uadapt}}$=\textcolor{red}{\zeta}$, {\ttfamily{isoadapt}}=5e-2, {\ttfamily{erradapt}}=1e-4, {\ttfamily{smoothadapt}}=5e-9, {\ttfamily{epsadapt}}=50e3.\\
\begin{figure}[!h]
\begin{center}
\begin{tabular}{m{7cm}m{7cm}}
    \includegraphics[height=7cm,width=7cm]{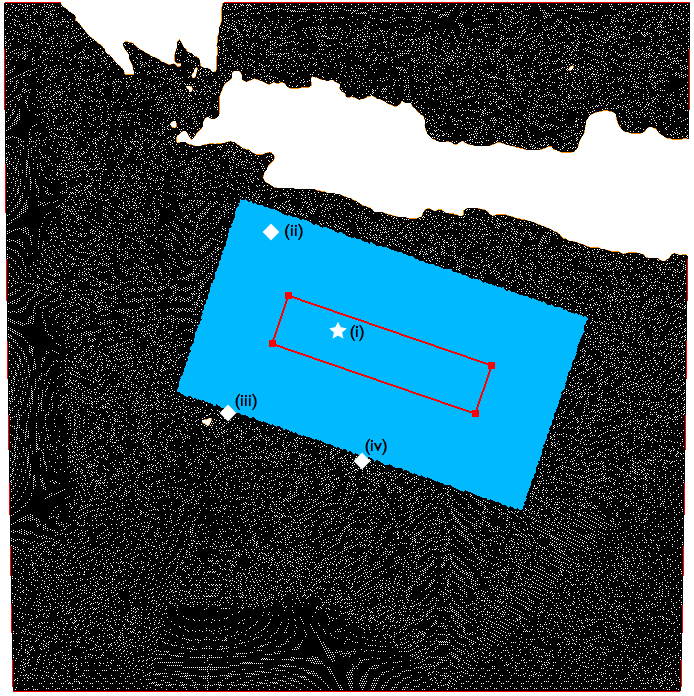} &\includegraphics[height=5cm,width=7cm]{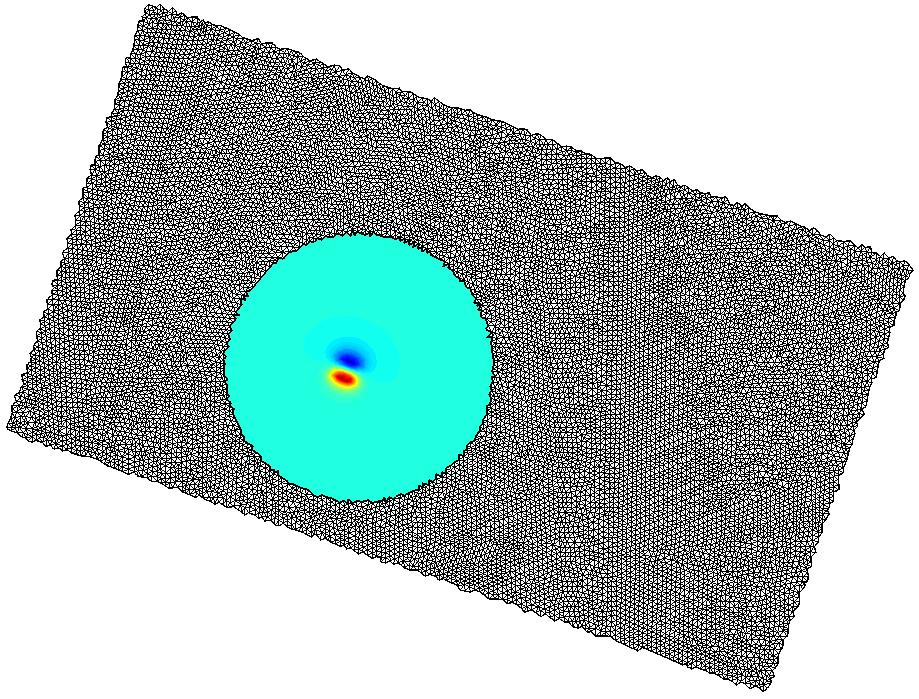}
\end{tabular}    
\caption{Left : Surface projection of the fault's plane and the mesh around, $\diamond$ : wave gauge, $\star$ : epicenter.\\
Right : the $14$-th Okada solution (min$=-0.09$ m, max$=0.17$ m ).\label{Fault_Okada_domain_mesh_sol} }
\end{center}
\end{figure}

We show in the Figure \ref{bottom_displ}, the bottom displacement $\textcolor{red}{\zeta}(x,y,t)$ at time $t=100$s and $t=270$s using our adapt mesh technique.

\begin{figure}[!h]
\begin{center}
\begin{tabular}{m{7cm}m{7cm}}
    \includegraphics[height=5cm,width=7cm]{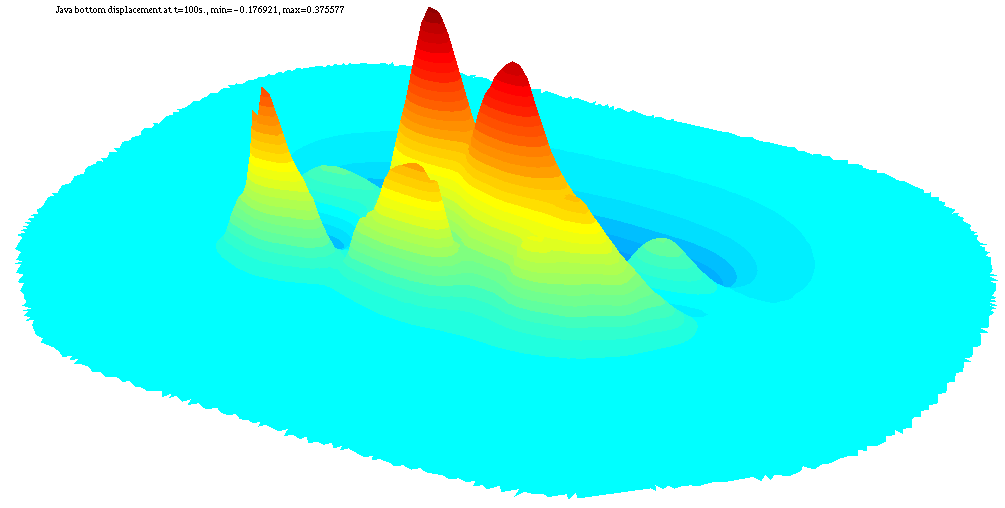} &\includegraphics[height=5cm,width=7cm]{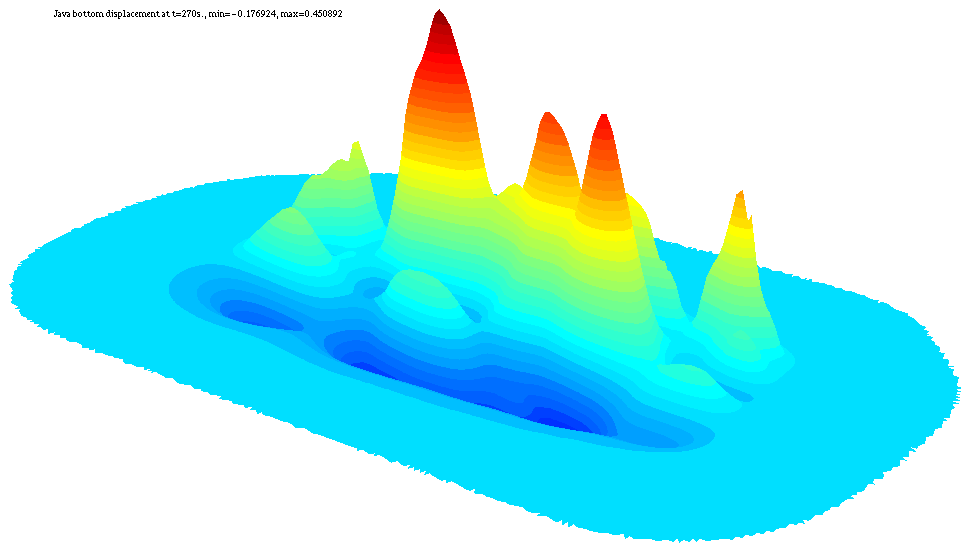}
\end{tabular}    
\caption{Bottom displacement at $t=100$ s (left, min$=-0.18$ m, max$=0.38$ m)  and at $t=270$ s (right, min$=-0.18$ m, max$=0.45$ m).\label{bottom_displ} }
\end{center}
\end{figure}
\begin{rem}
We note that after building the Okada solution $\mathcal{O}(x,y)$ in the {\it passive} generation or $\mathcal{O}_i(x,y)$ in the {\it Active} generation, we can remark that this solution is non-local and decays slowly to zero, that why in our adapt mesh technique we put 0 where the absolute value of the solution is less then $\min(\left|\min\left(\mathcal{O}_i(x,y)\right)\right|,\left|\max\left(\mathcal{O}_i(x,y)\right)\right|)/9.2$ meter, we make the same thing without adapt mesh in order to compare the solution using the same initial data.
\end{rem}
\section{Numerical simulations}\label{numsim}
In this section, we study first the rate of convergence of our codes for the sBBM (\ref{derivBBMFullsimp}) with non-dimensional and unscaled variable {\it i.e.}, with $g=1$ over a variable bottom in space, which establishes the adequacy of the chosen finite element discretization and the used time marching scheme, for the flat bottom case, we refer to \cite{Sad12a}\label{Sad12a0}, where we use the same technique as in this paper. Then we simulate the propagation of a wave, that looks like a {\it Tsunami} wave generated by an Earthquake, in the Mediterranean sea over the sBBM with a flat bottom, near Java island over a variable bottom in space and at the end near Java island over a variable bottom in space and in time.

In all numerical simulations we used $\mathbb{P}_1$ continuous piecewise linear functions for $\eta,u,v,\textcolor{blue}{D}$ and $\textcolor{red}{\zeta}.$
\subsection{Rate of convergence}
We prove in the figure below, that the RK2 time scheme considered for the sBBM variable bottom in space system is of order 2, we note that $\textcolor{red}{\zeta}(x,y,t)$ is only used for the generation of the {\it Tsunami} wave and will not be taken into account in the convergence rate test. In this example, we took Bi-Periodic Boundary Conditions for $\eta_{h}$, $ u_{h}$ and $ v_{h}$ on the whole boundary of the square $[0,2L]\times[0,2L]$, where $L=50$ and we consider the following exact solutions:
$$
\eta_{ex}=.2\cos(2\pi x/L-t)\cos(2\pi y/L-t),u_{ex}=.5\sin(2\pi x/L-t)\cos(2\pi y/L-t),
$$
$$
v_{ex}=.5\cos(2\pi x/L-t)\sin(2\pi y/L-t),\textcolor{blue}{D}(x,y)=1-.5\cos(2\pi x/L)\cos(2\pi y/L),
$$
adding an appropriate right-hand side function.\\

We measure at time $T=1$ and for $\theta^2=\ds\frac{2}{3},\ \delta t=\ds\frac{.01}{2^n}$ and $\delta x=\ds\frac{2L}{N}=\ds\frac{2L}{2^{n+5}}\ \forall n\in\{0,1,2,3,4\}$, the following errors 
$$
N_{L^2}(\eta)=\left\|\eta_h-\eta_{ex}\right\|_{L^2}, N_{H^1}(\eta)=\left\|\eta_h-\eta_{ex}\right\|_{H^1},N_{L^2}(\mathcal{V})=\left\|u_h-u_{ex}\right\|_{L^2}+\left\|v_h-v_{ex}\right\|_{L^2}
$$ 
$$
N_{H^1}(\mathcal{V})=\left\|u_h-u_{ex}\right\|_{H^1}+\left\|v_h-v_{ex}\right\|_{H^1}
$$
and we end up with the following results:\\

\begin{table}[!h]
\begin{center}
\begin{tabular}{|c|c|c|c|c|c|c|c|c|c|}
\hline N&$\delta t$ & $N_{L^2}(\eta)$ & rate & $N_{L^2}(\mathcal{V})$ & rate  & $N_{H^1}(\eta)$ & rate & $N_{H^1}(\mathcal{V})$ & rate% & time
\\ \hline $2^5$&$.01/2^0$ & 0.241446 & - & 1.10773 & - & 0.603174 & - & 1.62575 & -% & 136.852
\\ \hline $2^6$&$.01/2^1$ & 0.0607759 & 1.99013 & 0.280157 & 1.98329 & 0.301957 & 0.998228 & 0.812757 & 1.00021% & 1682.95
\\ \hline $2^7$&$.01/2^2$ & 0.01524 & 1.99564 & 0.0703759 & 1.99308 & 0.151186 & 0.998017 & 0.406962 & 0.99793% & 32475.5
\\ \hline $2^8$&$.01/2^3$ & 0.0038124 & 1.9987 & 0.017602 & 1.99909 & 0.075782 & 0.9975 & 0.203552 & 0.99965% & 72681.9
\\\hline
\end{tabular}
\end{center}
\caption{\label{rate_table} $L^2$ norm of the error for $\eta$ and $\mathcal{V}$.}
\end{table}
\begin{figure}[h!]
\begin{center}
\hspace{-1.25cm}\includegraphics[height=6cm]{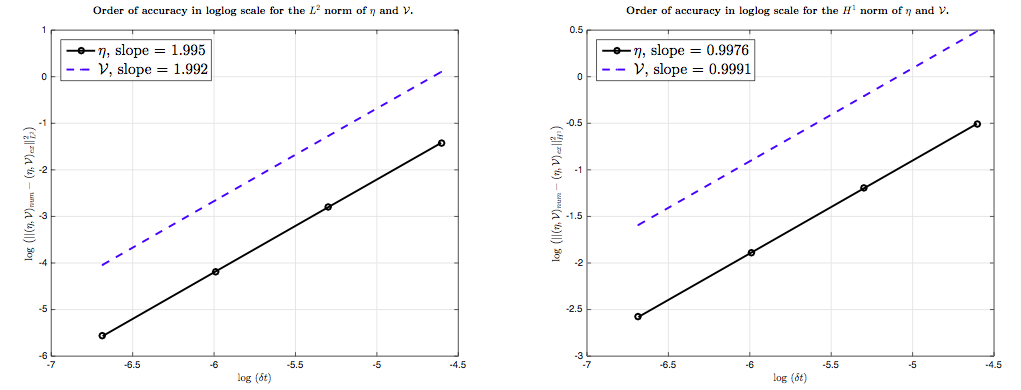}
\caption{\label{rate_sBBM} Rate of convergence for sBBM system with variable bottom in space.}
\end{center}
\end{figure}
So, the $L^2(\Omega \times ]0,T[)^2$ norm slope for $\eta$ and $\mathcal{V}$ is of order $\sim 2$ and the $L^2(0,T;H^1(\Omega)^2)$ slope for $\eta$ and $\mathcal{V}$ is of order $\sim 1$ as shown in the Figure \ref{rate_sBBM} and which confirms the convergence of the second-order Runge-Kutta scheme in time for the sBBM system with variable bottom in space.
\subsection{Propagation of a {\it Tsunami} wave in the Mediterranean sea with a flat bottom.}
In this section, the non-dimensional and unscaled variables in (\ref{derivBBMFullsimp}) {\it i.e.} $g=1$. We simulate here, the propagation of a wave that looks like a {\it Tsunami} wave generated by an Earthquake in the Mediterranean sea over the sBBM (\ref{derivBBMFullsimp}) with a flat bottom $-\textcolor{blue}{D}(x,y)=-1,5$ Km which is the average depth of the Mediterranean sea. This wave was defined above in the {\it passive} generation  part of the Section \ref{meshgeninitdata} where, in this case, the initial solution is concentrated in the small rectangle $[x_0-5W;x_0+4W]\times[y_0-1.5L;y_0+2.5L]$ and we take these following values : $L=20$ Km, $W=10$ Km, $\delta=7\degree$, $\phi=0\degree$, $\theta=90\degree$, $E=9,5$ GPA is the Young's modulus, $\nu=0,27$ is the Poisson's ratio, $U=2,5$ m, $(x_0;y_0)=(2390.*scale,590.*scale)$ and the fault depth$(x_0;y_0)=10$ Km. In this example, we will take the fact that the Lam\'e constants $\mu$ and $\lambda$ are given by the formulas $\mu=E/2(1+\nu)$ and $\lambda=E\nu/(1+\nu)(1-2\nu)$.\\
We also use the following settings : for the step time $\delta t=0.1=1$ s, a reflective BC for all the boundary, for the \textcolor{red}{\ttfamily{adaptmesh}} of \freefem :
\begin{lstlisting}[firstnumber=last]
fespace Vhinit(Thinit,P0);
Vhinit hT=hTriangle;
real Dx=hT[].min;
uadapt=eta0+u0+v0;
Th=adaptmesh(Th,uadapt,err=1.e-6,errg=1.e-2,hmin=Dx,iso=true,nbvx=1e8);
[eta0,u0,v0]=[eta0,u0,v0];	MAX=MAX;	D=D;
\end{lstlisting}
and for our new adapt mesh technique  :
\begin{lstlisting}[firstnumber=last]
fespace Wh(Th,P1);
mesh Thp=Th;
uadapt=eta0+u0+v0;
real isoadapt=5.e-2, erradapt=1.e-7,smoothadapt=5.e-3,epsadapt=2e-2;
bool output =true;
real wv=0.;
{ AdaptGS(Thp,uadapt,erradapt,isoadapt,smoothadapt,epsadapt,Thinit,output,wv); }
fespace Whp(Thp,P1);
Wh peta0=eta0,pu0=u0,pv0=v0;	Whp Dp,MAXp,eta0p,u0p,v0p;
matrix BWinterp=interpolate(Whp,Wh,inside=true);
eta0p[]=BWinterp*peta0[];u0p[]=BWinterp*pu0[];v0p[]=BWinterp*pv0[];MAXp[]=BWinterp*MAX[];
Th=Thp;
[eta0,u0,v0]=[eta0p,u0p,v0p];	MAX=MAXp;	D=D;
\end{lstlisting}
We note that, we adapt the mesh around the solution each 100 iterations {\it i.e.} each $10$ s by using the following value for the adapt mesh {\ttfamily{uadapt}}$=\eta+u+v$, {\ttfamily{isoadapt}}=5e-2, {\ttfamily{erradapt}}=1e-7, {\ttfamily{smoothadapt}}=5e-3, {\ttfamily{epsadapt}}=2e-2.
\begin{figure}[!h]
\begin{center}
\begin{tabular}{m{5cm}m{5cm}m{5cm}}
    \includegraphics[height=4.8cm,width=5cm]{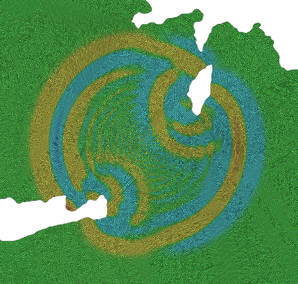} &\includegraphics[height=4.8cm,width=5cm]{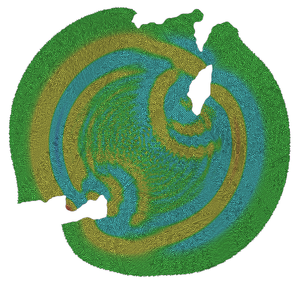}&\includegraphics[height=4.8cm,width=5cm]{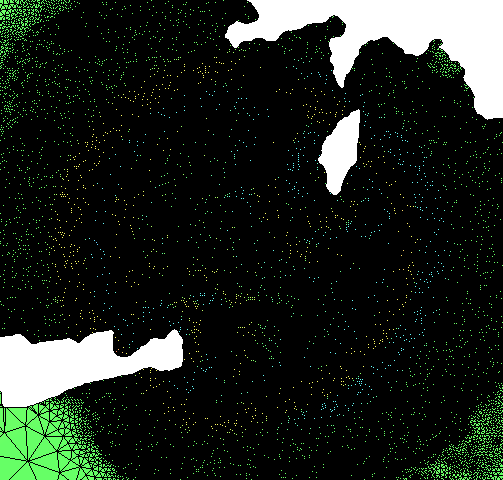}
\end{tabular}
\caption{The mesh and the solution at $t=1000$ s, with the Full method at left, the Adapt GS at the center and the Adapt FF with {\ttfamily{err=1.e-7}} at the right. \label{Tsu_Medit_prop_1} }
\end{center}
\end{figure}
\begin{figure}[!h]
\begin{center}
\begin{tabular}{m{5cm}m{5cm}m{5cm}}
    \includegraphics[height=4.8cm,width=5cm]{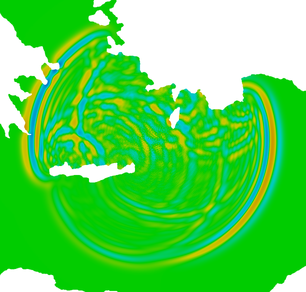} &\includegraphics[height=4.8cm,width=5cm]{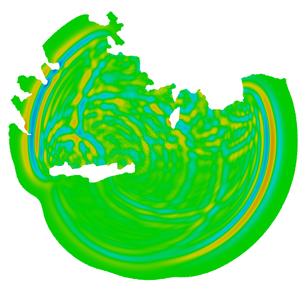}&\includegraphics[height=4.8cm,width=5cm]{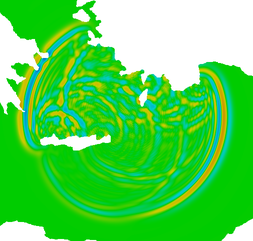}
\end{tabular}
\caption{The solution at $t=3000$ s, with the Full method at left, the Adapt GS at the center and the Adapt FF with {\ttfamily{err=1.e-7}} at the right. \label{Tsu_Medit_prop_2} }
\end{center}
\end{figure}

In order to compare the results between \textcolor{red}{\ttfamily{adaptmesh}} of \freefem, our new adapt mesh technique and without using mesh adaptation, we plot in addition to the free surface elevation $\eta$ in the Figures \ref{Tsu_Medit_prop_1} $\rightarrow$ \ref{Tsu_Medit_prop_2}, the variation of $\eta$ vs time in Figure \ref{Tsu_Medit_Gauge} at two wave 'gauges' placed at the positions represented by $\diamond$ in Figure \ref{Crete_mesh} (top, right) and the mass of the water $\int\eta$. Specifically, gauges were placed at the points $(i) : (2350.*scale,550.*scale)$, $(ii) : (2104.*scale,665.*scale)$.\\
In the Figure \ref{max_tsu_medit}, we represent the comparison between the three methods : Full, Adapt FF and Adapt GS of the maximum of the propagation of the solution at time $t=6800 sec$.\\
We also plot the computation time for each adapt mesh, the computation time of the simulation, the number of degree of freedom in Figure \ref{Tsu_Medit_comp}.\\
We can see in Figures \ref{Tsu_Medit_Gauge} and  \ref{Tsu_Medit_comp} that the \textcolor{red}{\ttfamily{adaptmesh}} of \freefem with {\ttfamily{err=1.e-2}} is the fastest method but unfortunately it does not preserve the mass invariant $\ds\int\eta$. On the other hand, our new adapt mesh technique preserves the mass invariant throughout the simulation with an error of order $2.1e-3$ and an important time computation difference with the one without mesh adaptation which is very good method for the {\it Tsunami} wave propagation.\\
For the \textcolor{red}{\ttfamily{adaptmesh}} of \freefem with {\ttfamily{err=1.e-7}}, we also get an almost a mass conservation with an error of order $9.5e-4$, but we obtain some difference in wave gauge with the Full method which is due to the refinement mesh adaptation and the interpolation of the solution, although the computation time is almost the double of the new adapt mesh technique.\\
We note also that we can go faster with our new mesh adaptation technique if we can also trunk the mass matrix after trunking the mesh, of course if the mass matrix is a constant along the simulation of the Full mesh, this is an outgoing project.
\begin{figure}[!h]
\begin{center}
    \hspace{-1cm}\includegraphics[height=8.5cm,width=18cm]{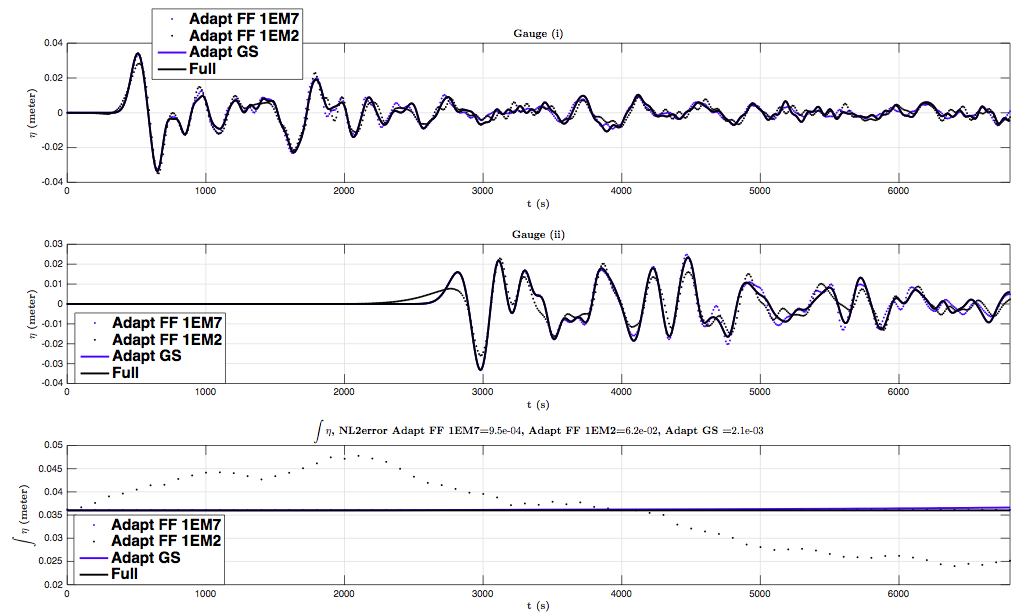}\\
\vspace{-1cm}
\caption{Comparison between the three methods : Full, Adapt FF and Adapt GS of the free surface elevations (in meters) vs time (in seconds), computed numerically at two wave gauges and of the mass conservation.\label{Tsu_Medit_Gauge} }
\end{center}
\end{figure}
\begin{figure}[h]
\begin{center}
    \hspace{-1cm}\includegraphics[height=4.5cm,width=18cm]{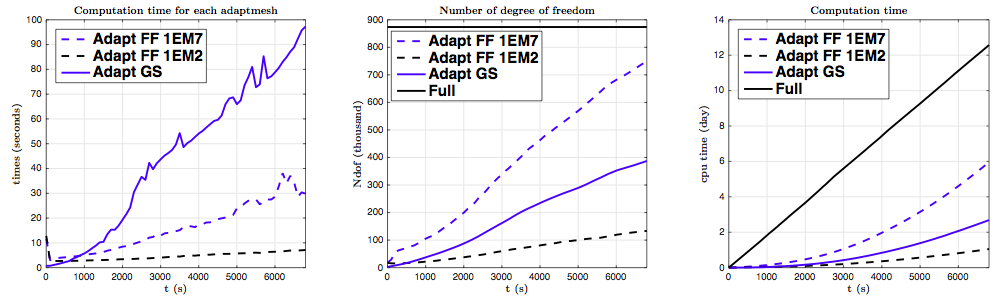}\\
\caption{Comparison between the three methods : Full, Adapt FF and Adapt GS of the computation time of each adaptemsh, the number of degree of freedom and the computation time of the simulation.\label{Tsu_Medit_comp} }
\end{center}
\end{figure}
\begin{figure}[!h]
\begin{center}
    \includegraphics[height=7cm,width=14.5cm]{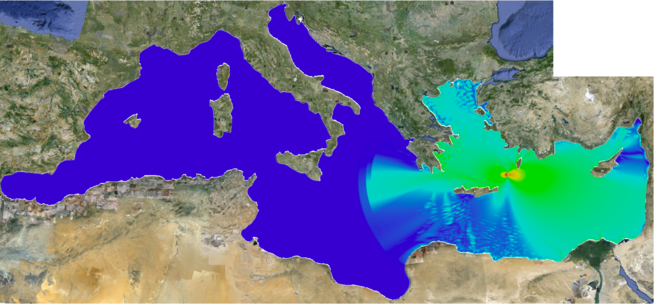}
        \includegraphics[height=7cm,width=14.5cm]{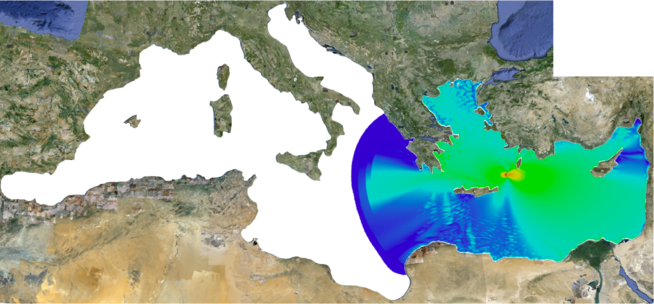}
\includegraphics[height=7cm,width=14.5cm]{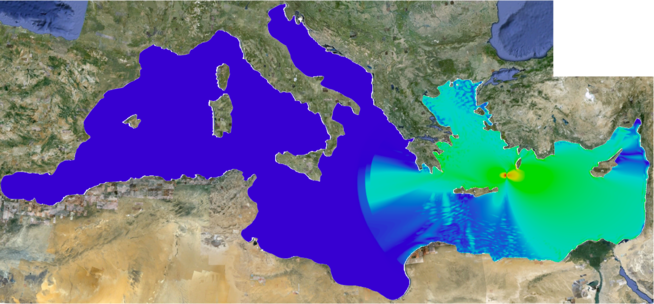}
\end{center}   
\caption{Comparison between the three method Full (up), Adapt GS (middle) and Adapt FF (down) of the maximum of the propagation of the solution of a {\it Tsunami} wave in the Mediterranean sea for $t=6800$ s.
\label{max_tsu_medit} 
}
\end{figure}
\clearpage
\subsection{Propagation of a {\it Tsunami} wave near Java island : {\it passive} generation .}
In this section, we will take the same initial data as defined above in the {\it passive} generation part of Section \ref{meshgeninitdata}, we take $\delta t=1$ s as the time step size and we note that, we adapt the mesh after computing the initial data for $\eta$ and then every 50 s by using the following value for the adapt mesh {\ttfamily{uadapt}}$=\eta+u+v$, {\ttfamily{isoadapt}}=3e-2, {\ttfamily{erradapt}}=1e-4, {\ttfamily{smoothadapt}}=5e-9, {\ttfamily{epsadapt}}=30e3.\\
We compare here the results between our new adapt mesh technique and without using mesh adaptation. To this end, we plot the free surface elevation $\eta$ in the Figures \ref{Tsu_Java_Pas_prop_1} and \ref{Tsu_Java_Pas_prop_2}.\\

\begin{figure}[!h]
\begin{center}
\begin{tabular}{m{5cm}m{5cm}m{5cm}}
    \includegraphics[height=5cm,width=5cm]{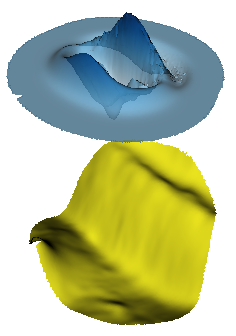} &\includegraphics[height=5cm,width=5cm]{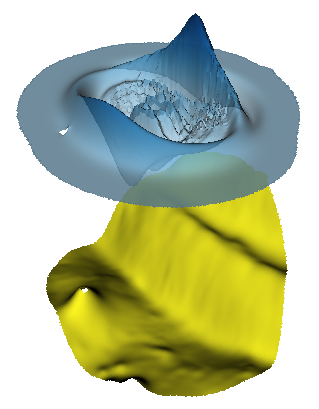}&\includegraphics[height=5cm,width=5cm]{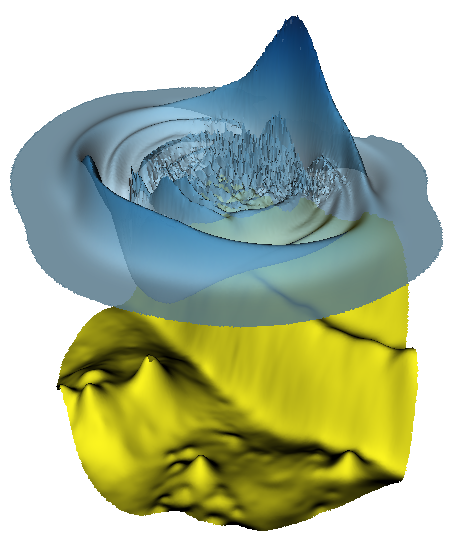}
\end{tabular}
\caption{{\it Passive} generation : The bottom and the solution at $t=250$ s (left, solution (min$=-0.36$ m, max$=0.38$ m), bottom (min$=-6207$ m, max$=-2096$ m)), $t=500$ s (center, solution (min$=-0.26$ m, max$=0.35$ m), bottom (min$=-6207$ m, max$=-243$ m)) and $t=1000$ s (right, solution (min$=-0.21$ m, max$=0.29$ m), bottom (min$=-6207$ m, max$=-100$ m)), with the Adapt GS method. \label{Tsu_Java_Pas_prop_1} }
\end{center}
\end{figure}
\begin{figure}[!h]
\begin{center}
\begin{tabular}{m{8cm}m{8cm}}
    \includegraphics[height=7cm,width=8cm]{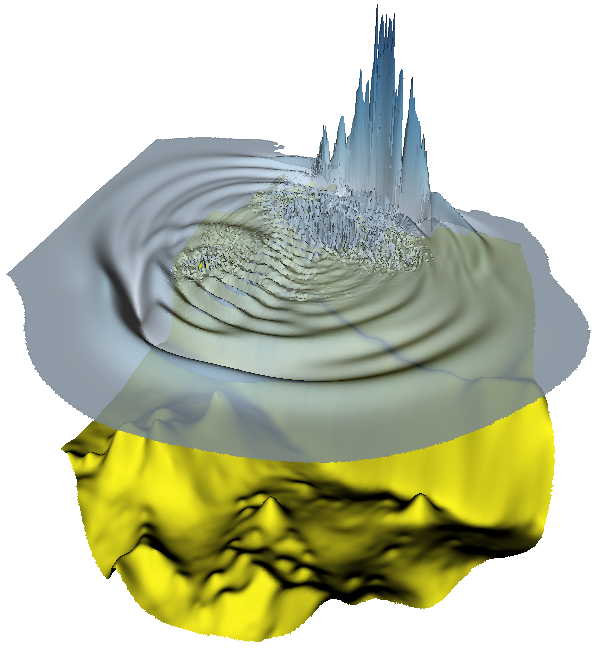} &\includegraphics[height=7cm,width=8cm]{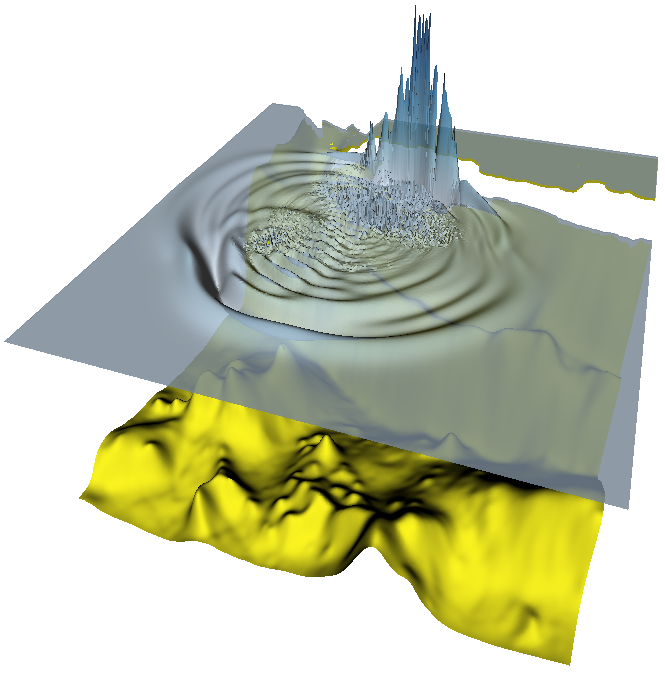}
\end{tabular}
\caption{{\it Passive} generation : Comparison between the bottom and the solution at $t=1500$ s, with the Adapt GS method (left, solution (min$=-0.38$ m, max$=1.07$ m), bottom (min$=-6207$ m, max$=-100$ m)) and with the Full one (right, solution (min$=-0.38$ m, max$=1.07$ m), bottom (min$=-6207$ m, max$=-100$ m)). \label{Tsu_Java_Pas_prop_2} }
\end{center}
\end{figure}
\vspace{-.5cm}
We also plot, the variation of $\eta$ vs time (in Figure \ref{gauge_Java_Passive}) at four numerical wave gauges placed at the following locations: (i) $(107.345\degree,-9.295\degree)$, (ii) $(106.5\degree,-8\degree)$, (iii) $(105.9\degree,-10.35\degree)$ and (iv) $(107.7\degree,-11\degree)$ (see Figure \ref{Fault_Okada_domain_mesh_sol} (left)) where (i) is the position of the epicenter. However, because of the large variations of the bottom, shorter waves were generated, especially around Christmas Island (southwest of Java) and around the undersea canyon near the Earthquake's epicenter.

\begin{figure}[h]
\begin{center}
    \hspace{-1cm}\includegraphics[height=7.8cm,width=18cm]{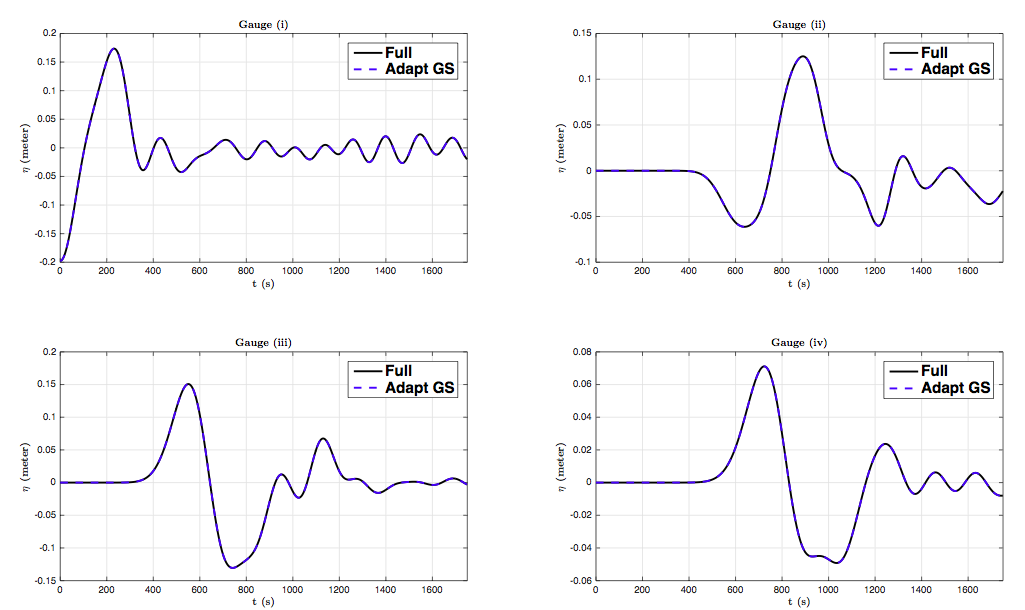}\\
\caption{{\it passive} generation : Comparison between the two methods Full and Adapt GS of the free surface elevations (in meters) vs time (in seconds), computed numerically at four wave gauges where the gauge (i) correspond to the epicenter.\label{gauge_Java_Passive} }
\end{center}
\end{figure}
Finally, we present a comparison of the Kinetic, Potential and Total energy with the Full mesh (in Figure \ref{comp_Java_Passive}, top left) and with the Adapt GS method (in Figure \ref{comp_Java_Passive}, top right) defined in \cite{DD09}\label{DD090} as:
$$
E_c=\frac{1}{2}\rho_w\int_{\Omega}\left(\int_{-D(x,y)}^{\eta}|V|^2\mbox{d}z\right)\mbox{d}x\mbox{d}y, \qquad E_p=\frac{1}{2}\rho_w\cdot g\int_{\Omega}\eta^2\mbox{d}x\mbox{d}y,
$$
where $\rho_w=1027$ $\mbox{Kg/m}^3$ is the ocean water density, the number of degree of freedom (in Figure \ref{comp_Java_Passive}, down left) and the computation time of the simulation (in Figure \ref{comp_Java_Passive}, down right). We obtain here an error of order $2.6e-4$ between the Total Energy with Adapt GS and without adaptation.
\begin{figure}[h]
\begin{center}
    \hspace{-1cm}\includegraphics[height=7.8cm,width=18cm]{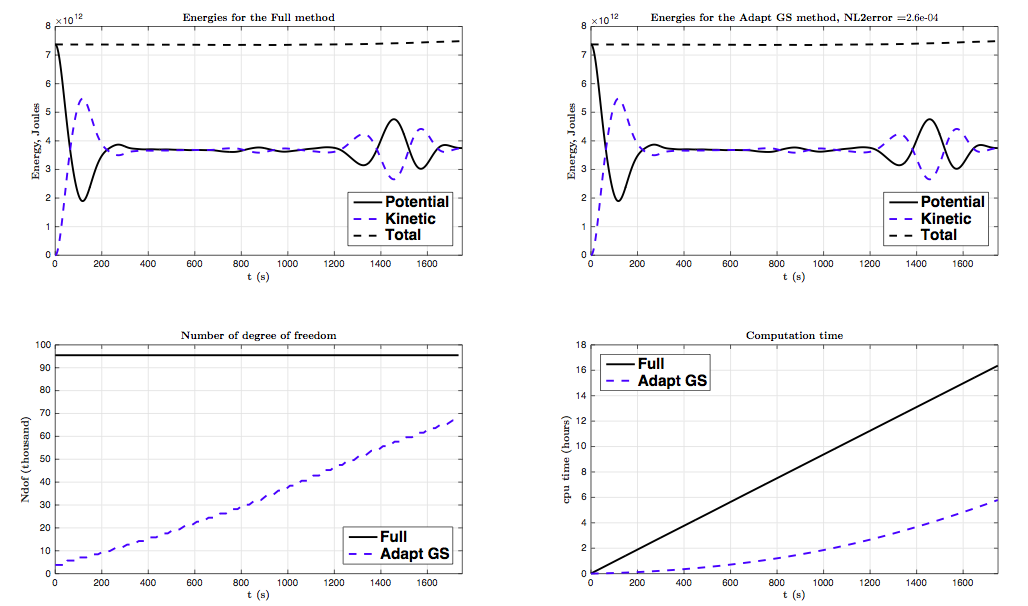}\\
\caption{{\it passive} generation : Comparison between the two methods Full and Adapt GS of the Kinetic, Potential and Total energy, the number of degree of freedom and the computation time of the simulation.\label{comp_Java_Passive} }
\end{center}
\end{figure}

We present in the Figure \ref{max_Tsu_Java_Passive} the comparison of the maximum of the propagation of the solution between the Full and the Adapt GS method at $t=1750$ s.
\clearpage
\begin{figure}[!h]
\begin{center}
\begin{tabular}{m{8cm}m{8cm}}
    \includegraphics[height=8cm,width=8cm]{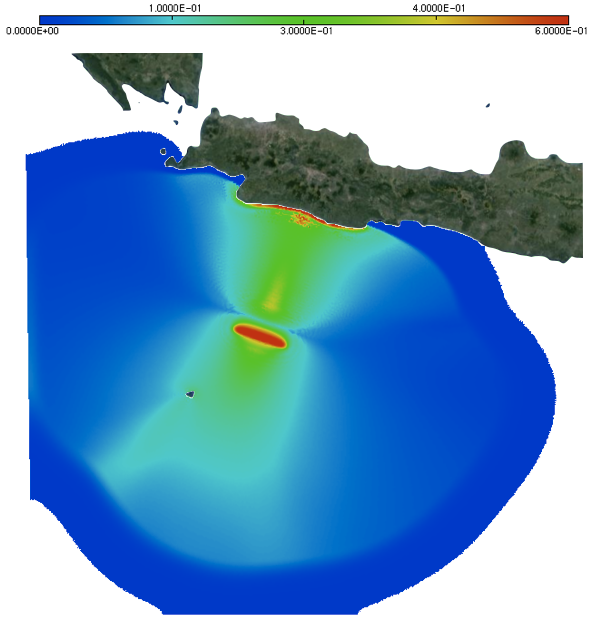} &\includegraphics[height=8cm,width=8cm]{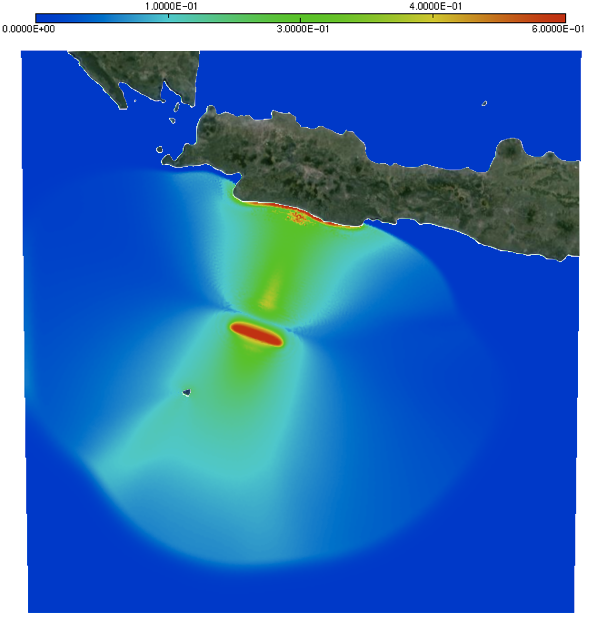}
\end{tabular}
\caption{{\it Passive} generation : Comparison between the maximum of the solution at $t=1750$ s, with the Adapt GS method (left) and with the Full one (right). \label{max_Tsu_Java_Passive} }
\end{center}
\end{figure}
\subsection{Propagation of a {\it Tsunami} wave near Java island : {\it active} generation .}
For a more realistic case as in the Java 2006 event, we use the {\it active} generation in order to model the generation of a {\it Tsunami} wave as in \cite{DMGD13,DMCS12}\label{DMGD134}\label{DMCS120}. In this case we consider zero initial conditions for both the surface elevation and the velocity field, we take $\delta t=2$ s as the time step size, we assume that the bottom described in the Section \ref{meshgeninitdata} is moving in time and we note that we adapt the mesh, before the end of the generation time $t=270$ s, each 3 iterations {\it i.e.} each $6$ s by using the following value for the adapt mesh {\ttfamily{uadapt}}$=\eta+u+v$, {\ttfamily{isoadapt}}=5e-2, {\ttfamily{erradapt}}=1e-4, {\ttfamily{smoothadapt}}=5e-9, {\ttfamily{epsadapt}}=50e3 and then for $t>270$ s each 25 iterations {\it i.e.} each $50$ s.\\
We compare here only the results between our new adapt mesh technique and without using mesh adaptation. To this end, we plot the free surface elevation $\eta$ in the Figures \ref{Tsu_Java_Act_prop_1} $\rightarrow$ \ref{Tsu_Java_Act_prop_3}. However, because of the large variations of the bottom, shorter waves were generated, especially around Christmas Island (southwest of Java) and around the undersea canyon near the Earthquake's epicenter.\\
We also plot, the variation of $\eta$ vs time (in Figure \ref{gauge_Java_Active}) at four numerical wave gauges placed at the following locations: (i) $(107.345\degree,-9.295\degree)$, (ii) $(106.5\degree,-8\degree)$, (iii) $(105.9\degree,-10.35\degree)$ and (iv) $(107.7\degree,-11\degree)$ (see Figure \ref{Fault_Okada_domain_mesh_sol} (left)) where (i) is the position of the epicenter.\\
\begin{figure}[!h]
\begin{center}
\begin{tabular}{m{5cm}m{5cm}m{5cm}}
    \includegraphics[height=5cm,width=5cm]{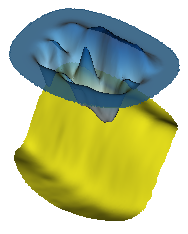} &\includegraphics[height=5cm,width=5cm]{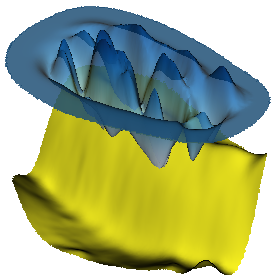}&\includegraphics[height=5cm,width=5cm]{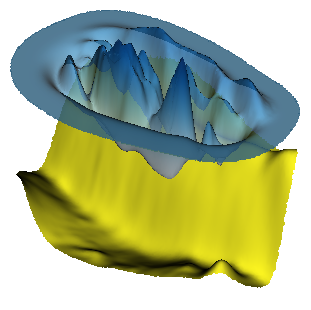}
\end{tabular}
\caption{{\it Active} generation : The bottom and the solution at $t=100$ s (left, solution (min$=-0.17$ m, max$=0.07$ m), bottom (min$=-6207$ m, max$=-2589$ m)), $t=200$ s (center, solution (min$=-0.19$ m, max$=0.08$ m), bottom (min$=-6207$ m, max$=-2285$ m)) and $t=270$ s (right, solution (min$=-0.14$ m, max$=0.10$ m), bottom (min$=-6207$ m, max$=-2084$ m)), with the Adapt GS method. \label{Tsu_Java_Act_prop_1} }
\end{center}
\end{figure}
\begin{figure}[!h]
\begin{center}
\begin{tabular}{m{8cm}m{8cm}}
    \includegraphics[height=6cm,width=7cm]{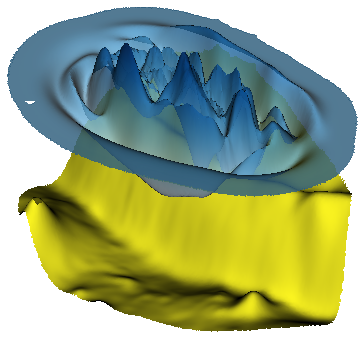} &\includegraphics[height=6cm,width=7cm]{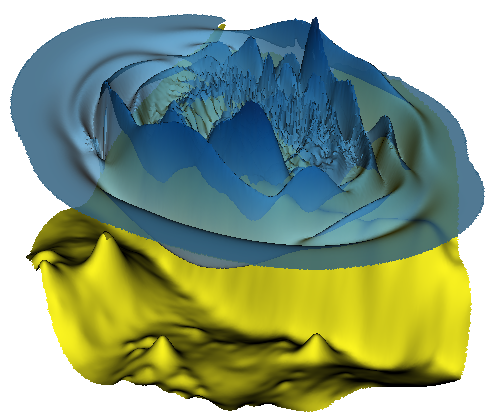}
\end{tabular}
\caption{{\it Active} generation : The bottom and the solution at $t=500$ s (left, solution (min$=-0.15$ m, max$=0.10$ m), bottom (min$=-6207$ m, max$=-260$ m)) and $t=1000$ s (right, solution (min$=-0.14$ m, max$=0.09$ m), bottom (min$=-6207$ m, max$=-100$ m)), with the Adapt GS method. \label{Tsu_Java_Act_prop_2} }
\end{center}
\end{figure}

\begin{figure}[!h]
\begin{center}
\begin{tabular}{m{8cm}m{8cm}}
    \includegraphics[height=7cm,width=8cm]{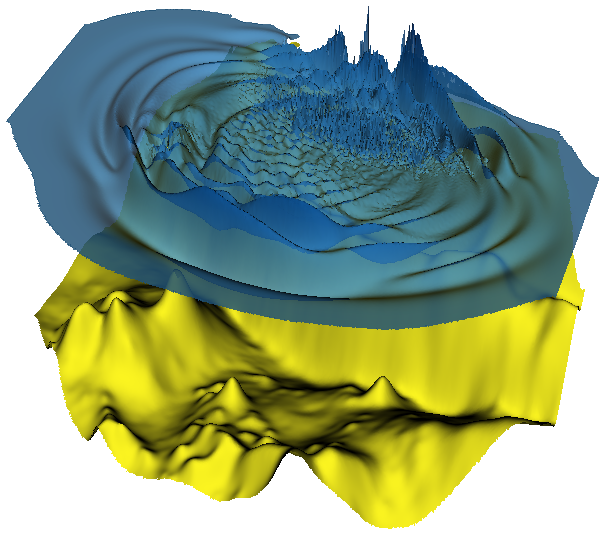} &\includegraphics[height=7cm,width=8cm]{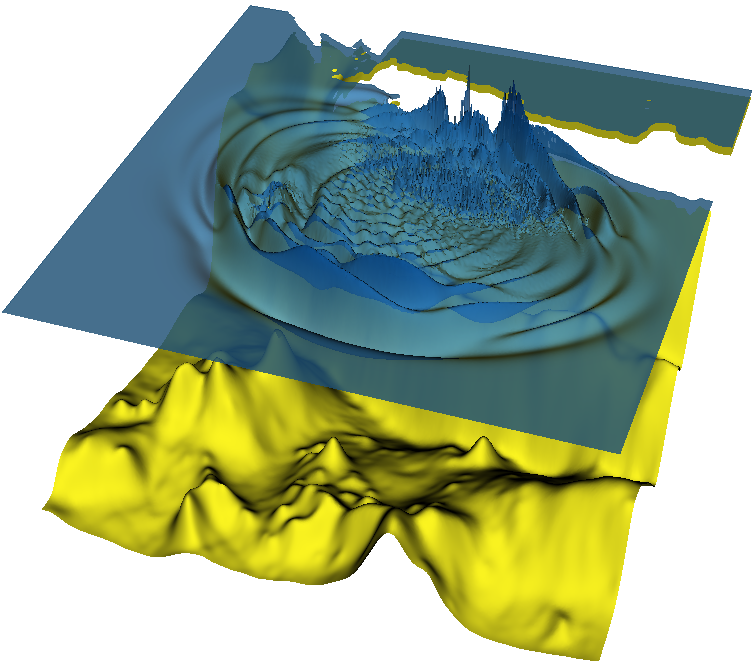}
\end{tabular}
\caption{{\it Active} generation : Comparison between the bottom and the solution at $t=1500$ s, with the Adapt GS method (left, solution (min$=-0.29$ m, max$=0.13$ m), bottom (min$=-6207$ m, max$=-100$ m)) and with the Full one (right, solution (min$=-0.29$ m, max$=0.13$ m), bottom (min$=-6207$ m, max$=-100$ m)). \label{Tsu_Java_Act_prop_3} }
\end{center}
\end{figure}
At the end, we present a comparison of the Kinetic, Potential and Total energy with the Full mesh (in Figure \ref{comp_Java_Active}, top left) and with the Adapt GS method (in Figure \ref{comp_Java_Active}, top right) defined in \cite{DD09}\label{DD090} as:
$$
E_c=\frac{1}{2}\rho_w\int_{\Omega}\left(\int_{-D(x,y)}^{\eta}|V|^2\mbox{d}z\right)\mbox{d}x\mbox{d}y, \qquad E_p=\frac{1}{2}\rho_w\cdot g\int_{\Omega}\eta^2\mbox{d}x\mbox{d}y,
$$
the number of degree of freedom (in Figure \ref{comp_Java_Active}, down left) and the computation time of the simulation (in Figure \ref{comp_Java_Active}, down right). We obtain here an error of order $2e-5$ between the Total Energy with Adapt GS and without adaptation.\\
We present in the Figure \ref{max_Tsu_Java_Active} the comparison of the maximum of the propagation of the solution between the Full and the Adapt GS method at $t=1750$ s.
\clearpage
\begin{figure}[h]
\begin{center}
    \hspace{-1cm}\includegraphics[height=8cm,width=18cm]{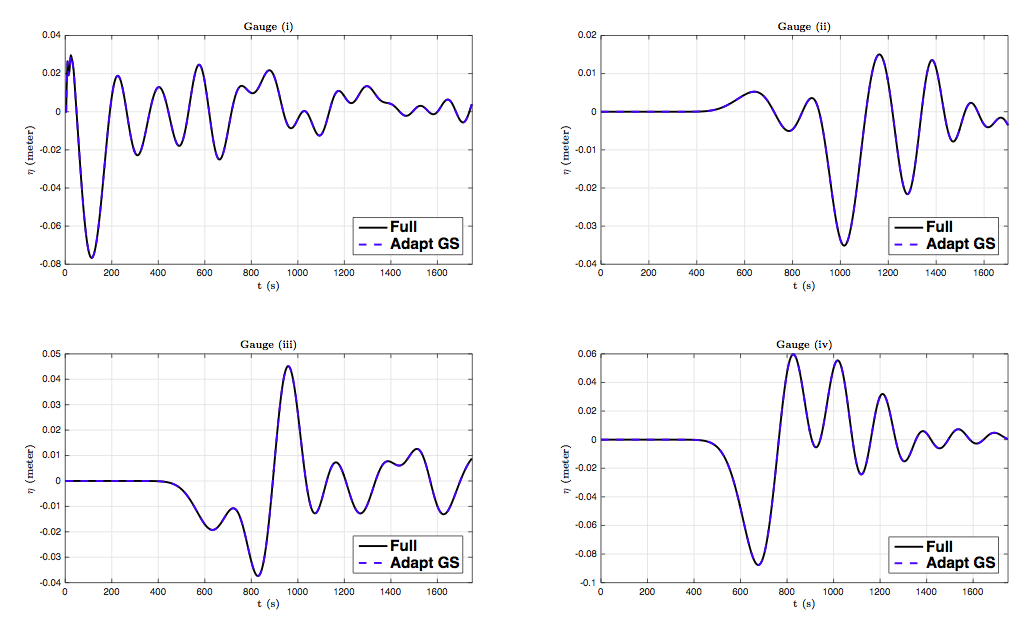}\\
\caption{{\it Active} generation : Comparison between the two methods Full and Adapt GS of the free surface elevations (in meters) vs time (in seconds), computed numerically at four wave gauges where the gauge (i) correspond to the epicenter.\label{gauge_Java_Active} }
\end{center}
\end{figure}
\begin{figure}[!h]
\begin{center}
    \hspace{-1cm}\includegraphics[height=10cm,width=18cm]{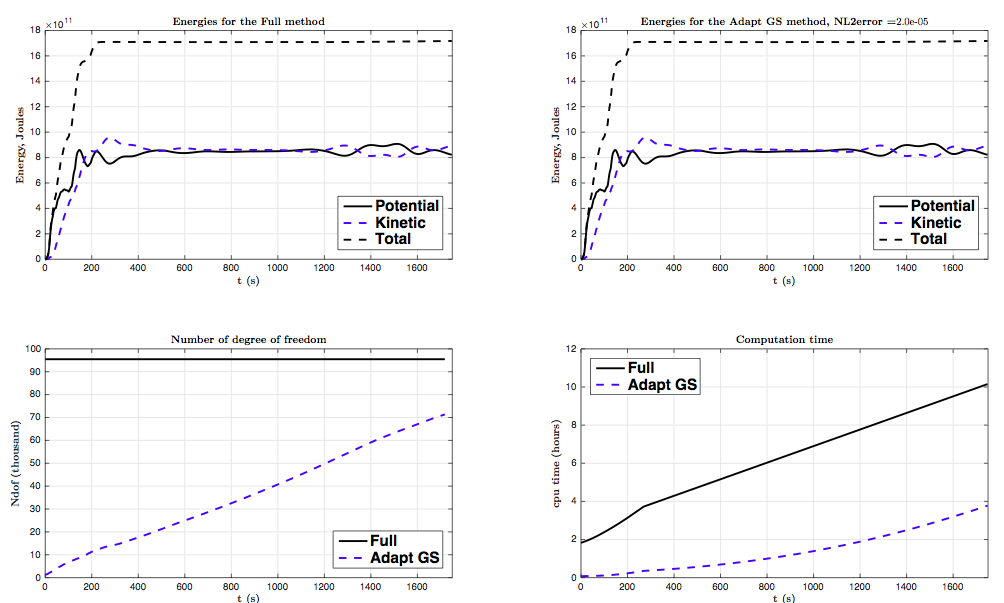}\\
\caption{{\it Active} generation : Comparison between the two methods Full and Adapt GS of the Kinetic, Potential and Total energy, the number of degree of freedom and the computation time of the simulation.\label{comp_Java_Active} }
\end{center}
\end{figure}
\clearpage
\begin{figure}[!h]
\begin{center}
\begin{tabular}{m{8cm}m{8cm}}
    \includegraphics[height=8cm,width=8cm]{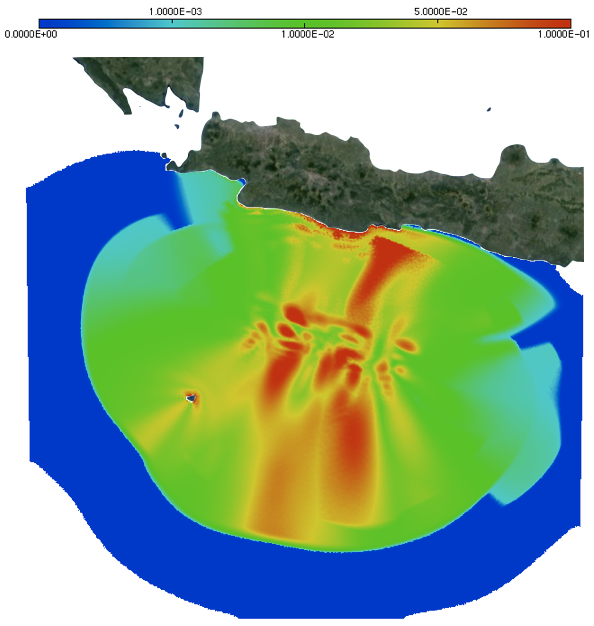} &\includegraphics[height=8cm,width=8cm]{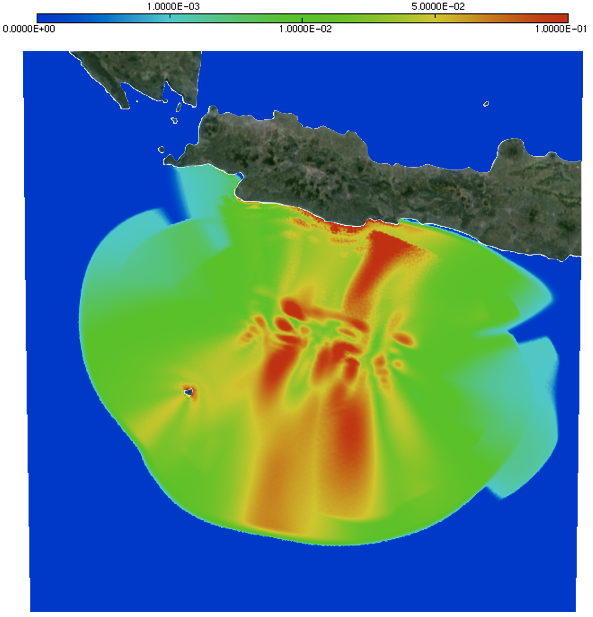}
\end{tabular}
\caption{{\it Active} generation : Comparison between the maximum of the solution at $t=1750$ s, with the Adapt GS method (left) and with the Full one (right). \label{max_Tsu_Java_Active} }
\end{center}
\end{figure}
\section{Conclusion and Outlook}
We show in this paper, the usefulness of \freefem for the {\it simplified} Boussinesq system of BBM type by building the domain, on the one hand using a photo taken from Google Earth$\circledR$ and on the other hand through an xyz bathymetry downloaded from the NOAA website. For the simulation of a {\it Tsunami} wave near Java island, the digital computing environment that we developed allows the integration of realistic data (bathymetry and geography) in a relatively simple framework. Another work concerning the generation, propagation and inundation of a {\it Tsunami} wave will be discussed in the case of the Shallow Water equations in \cite{Sad12b}\label{Sad12b0}, where in this case, we are not constraint by the smoothness of the bathymetry to avoid blow-up and where the same special adapt technique with a parallel version of the code will be introduced.\\

{\bf Acknowledgements : } This work would not be done without a lot of advice and special care from Denys Dutykh, Fr\'ed\'eric Hecht, Dimitrios Mitsotakis, and Olivier Pantz, to whom I am grateful for all their help, their fruitful discussions and their remarks.\\

All the videos of the simulations of a {\it Tsunami} wave for the results presents in this paper are given in the following links :\\
\url{http://www.lamfa.u-picardie.fr/sadaka/movies/Tsu_Medit_Full.mov}\\
\url{http://www.lamfa.u-picardie.fr/sadaka/movies/Tsu_Medit_Adapt_GS.mov}\\
\url{http://www.lamfa.u-picardie.fr/sadaka/movies/Tsu_Medit_Adapt_FF_1EM2_sol.mov}\\
\url{http://www.lamfa.u-picardie.fr/sadaka/movies/Tsu_Medit_Adapt_FF_1EM2_mesh.mov}\\
\url{http://www.lamfa.u-picardie.fr/sadaka/movies/Tsu_Medit_Adapt_FF_1EM7_sol.mov}\\
\url{http://www.lamfa.u-picardie.fr/sadaka/movies/Tsu_Medit_Adapt_FF_1EM7_mesh.mov}\\
\url{http://www.lamfa.u-picardie.fr/sadaka/movies/Okada_Java_Dynamic.gif}\\
\url{http://www.lamfa.u-picardie.fr/sadaka/movies/Java_Bottom_Displacement.gif}\\
\url{http://www.lamfa.u-picardie.fr/sadaka/movies/Tsu_Java_sBBM_Pas_Adapt_GS_Full.mov}\\
\url{http://www.lamfa.u-picardie.fr/sadaka/movies/Tsu_Java_sBBM_Act_Adapt_GS_Full.mov}

\bibliographystyle{plain}

}
\end{document}